\documentstyle[amstex]{amsart}

\makeatletter

\def\cicy#1(#2|#3)#4{\left(\matrix{#2}\right|\!\!
		     \left|\matrix{#3}\right)^{{#4}}_{#1}}

\begin{document}

\newcommand{\e}{\epsilon}
\newcommand{\w}{{\bold w}}
\newcommand{\y}{{\bold y}}
\newcommand{\z}{{\bold z}}
\newcommand{\x}{{\bold x}}
\newcommand{\N}{{\bold N}}
\newcommand{\Z}{{\bold Z}}
\newcommand{\F}{{\bold F}}
\newcommand{\R}{{\bold R}}
\newcommand{\Q}{{\bold Q}}
\newcommand{\C}{{\bold C}}
\newcommand{\BP}{{\bold P}}
\newcommand{\cO}{{\mathcal O}}
\newcommand{\cX}{{\mathcal X}}
\newcommand{\cH}{{\mathcal H}}
\newcommand{\cM}{{\mathcal M}}
\newcommand{\sB}{{\sf B}}
\newcommand{\cT}{{\mathcal T}}
\newcommand{\cI}{{\mathcal I}}
\newcommand{\cS}{{\mathcal S}}
\newcommand{\sE}{{\sf E}}

\newcommand{\sA}{{\sf A}}
\newcommand{\ga}{{\sf a}}
\newcommand{\es}{{\sf s}}
\newcommand{\m}{{\bold m}}
\newcommand{\bS}{{\bold S}}

\newcommand{\ihra}{\stackrel{i}{\hookrightarrow}}
\newcommand\rank{\mathop{\rm rank}\nolimits}
\newcommand\im{\mathop{\rm Im}\nolimits}
\newcommand\coker{\mathop{\rm coker}\nolimits}
\newcommand\Li{\mathop{\rm Li}\nolimits}
\newcommand\NS{\mathop{\rm NS}\nolimits}
\newcommand\Hom{\mathop{\rm Hom}\nolimits}
\newcommand\Ext{\mathop{\rm Ext}\nolimits}
\newcommand\Pic{\mathop{\rm Pic}\nolimits}
\newcommand\Spec{\mathop{\rm Spec}\nolimits}
\newcommand\Hilb{\mathop{\rm Hilb}\nolimits}
\newcommand{\length}{\mathop{\rm length}\nolimits}

\newcommand\lra{\longrightarrow}
\newcommand\ra{\rightarrow}
\newcommand\la{\leftarrow}
\newcommand\cJ{{\mathcal J}}
\newcommand\JG{J_{\Gamma}}
\newcommand{\wvskp}{\vspace{1cm}}
\newcommand{\vskp}{\vspace{5mm}}
\newcommand{\nvskp}{\vspace{1mm}}
\newcommand{\nid}{\noindent}
\newcommand{\new}{\nvskp \nid}
\newtheorem{Assumption}{Assumption}[section]
\newtheorem{Theorem}{Theorem}[section]
\newtheorem{Lemma}{Lemma}[section]
\newtheorem{Remark}{Remark}[section]
\newtheorem{Corollary}{Corollary}[section]
\newtheorem{Conjecture}{Conjecture}[section]
\newtheorem{Proposition}{Proposition}[section]
\newtheorem{Example}{Example}[section]
\newtheorem{Definition}{Definition}[section]
\newtheorem{Question}{Question}[section]
\renewcommand{\thesubsection}{\it}

\baselineskip=14pt

\title{Local cohomology of generalized Okamoto--Painlev\'e pairs and Painlev\'e equations}
\author{Hitomi Terajima}
\address{Department of Mathematics, Faculty of Science, 
Kobe University, Kobe, Rokko, 657-8501, Japan}
\email{terajima@math.kobe-u.ac.jp}
\maketitle

\begin{abstract}
In the theory of deformation of Okamoto-Painlev\'e pair $(S,Y)$, a local cohomology group $H^1_D(\Theta_S(-\log D))$ plays an important role.
In this paper, we estimate the local cohomology group of pair $(S,Y)$ for several types, and obtain the following results. For a pair $(S,Y)$ corresponding to the space of initial conditions of the Painlev\'e equations, we show that the local cohomology group $H^1_D(\Theta_S(-\log D))$ is at least 1 dimensional. This fact is the key to understand Painlev\'e equation related to $(S,Y)$. Moreover we show that, for the pairs $(S,Y)$ of type $\tilde{A_8}$, the local cohomology group $H^1_D(\Theta_S(-\log D))$ vanish. Therefore in this case, there is no differential equation on $S-Y$ in the sense of the theory.
\end{abstract}

\section{Introduction}

To study the spaces of initial conditions of Painlev\'e equations constructed by Okamoto \cite{O1} \cite{O2} \cite{O3}, we introduced the notion of {\it generalized Okamoto-Painlev\'e pair} $(S,Y)$ in \cite{STT}. This is a pair of a complex projective surface $S$ and an anti-canonical divisor $Y \in |-K_S|$ of $S$ satisfying the following conditions: For the irreducible decomposition $Y=\sum^r_{i=1} m_i Y_i$, one has $Y \cdot Y_i = \deg Y|_{Y_i} = 0$ for $1 \leq i \leq r$. In addition, if $S$ is a rational surface, $(S,Y)$ is called a {\it generalized rational Okamoto-Painlev\'e pair}.

The generalized rational Okamoto-Painlev\'e pairs of non-fibered type are classified into three types; elliptic, multiplicative and additive. Moreover, the generalized rational Okamoto-Painlev\'e pairs of additive type such that $D:=Y_{red}=\sum^r_{i=1}Y_i$ is a divisor with only normal crossings correspond to the spaces of initial conditions of the Painlev\'e equations (cf. \cite{STT} \cite{Sakai}) as follows:\vspace{2mm}
\begin{center}
\renewcommand{\arraystretch}{1.5}
\begin{tabular}{|c||c|c|c|c|c|c|c|c|} 
\hline
	type of $(S,Y)$ & $\tilde{E_8}$ & $\tilde{E_7}$ & $\tilde{D_8}$ & $\tilde{D_7}$ & 
$\tilde{D_6}$ & $\tilde{E_6}$ &$\tilde{D_5}$ & $\tilde{D_4}$ \\
	\hline
	Painlev\'{e} equations & $P_{I}$  & $P_{II}$ & $P^{\tilde{D_8}}_{III}$ & $P^{\tilde{D_7}}_{III}$ & $P_{III}=P^{\tilde{D_6}}_{III}$ & $P_{IV}$ & $P_{V}$  &  $P_{VI}$ \\
	\hline
\end{tabular}
\end{center}\vspace{2mm}

In what follows,  $(S,Y)$ is a generalized rational Okamoto-Painlev\'e pair of non-fibered type satisfying the condition: $D=Y_{red}$ is a normal crossing divisor with at least two irreducible components so that all irreducible components of $Y_{red}$ are smooth rational curves.  Let $\Theta_S( - \log D )$ be the sheaf of regular vector fields which have logarithmic zero along $D$. Here and after, all sheaves of $\cO_S$-modules are considered in algebraic category.

Then we have the following key exact sequence:
$$
0 \ra  H^1_D(\Theta_S( - \log D))  \ra  H^1(S, \Theta_S( - \log D)) \stackrel{res}{\ra} H^1(S - D, \Theta_S( - \log D)) ,
$$
where $H^1(S, \Theta_S( - \log D))$ and $H^1(S - D, \Theta_S( - \log D))$ are the space of infinitesimal deformations of the pair $(S,Y)$ (cf. \cite{Kaw}) and the space of infinitesimal deformations of $S-Y$, respectively.\\
In \cite{STT}, we show that the directions corresponding to local cohomology $H^1_D(\Theta_S(-\log D))$ in the deformation of a pair $(S,Y)$ induce differential equations on $S-Y$, by generalizing the Kodaira-Spencer theory to the open surface $S-Y$.

In this paper, we will show 
$$
\dim H^1_D(\Theta_S(-\log D)) \geq 1,
$$
when $(S,Y)$ is of additive type with the normal crossing divisor $D=Y_{red}$ (Theorem \ref{thm:time}). This result is natural since $S-Y$ corresponds to the space of initial conditions of the Painlev\'e equations in this case.
On the other hand, this is not always the case. In fact, we will prove that
$$
H^1_D(\Theta_S(-\log D))=\{ 0 \},
$$
for pairs $(S,Y)$ of type $\tilde{A_8}$ which is classified as a multiplicative type (Proposition \ref{prop:multi}). This means that there does not exist differential equation on $S-Y$ by the way above. All these computations are carried out through calculations of \v{C}ech cohomologies by taking a coordinate system explicitly.\\

\section{Local cohomology sequences and Time variables}\label{sec:local}

Let $(S, Y)$ be a generalized rational  Okamoto--Painlev\'{e} pair and set $D = Y_{red}$. Moreover, in this section, 
we assume that
\begin{enumerate}

\item $(S, Y)$ is  of non-fibered type and 

\item  $Y_{red}$ is a simple normal crossing divisor with 
at least two irreducible components,  i.e. $(r \geq 2)$ 
so that all irreducible components of $Y_{red}$ are 
smooth rational curves.  
\end{enumerate} 

Here $(S,Y)$ is called of fibered type if $S$ has a structure of an elliptic surfacefibration $f:S \to \BP^1$ with $f^{\ast} (\infty )=n\, Y$ for some $n \geq 1$. If $(S,Y)$ is not of fibered type, we call $(S,Y)$ of non-fibered type. (cf. Definition 1.3, \cite{STT}).

In what follows, $\cO_S$ and  $\cO_{S-D}$ denote the sheaves
of germs of algebraic regular functions on $S$ and $S - D$ 
respectively. Moreover all sheaves of $\cO_S$-modules are 
considered in algebraic category unless otherwise stated.   
Let us consider the following exact sequence of local 
cohomology groups ([Corollary 1.9, [Gr]])
$$
\begin{array}{cccccc}
H^0(S, \Theta_S(- \log D )) &  \ra & H^0(S-D, \Theta_{S}( - \log D)) & \ra & H^1_D(\Theta_S( - \log D)) & \ra \\ 
 H^1(S, \Theta_S( - \log D)) & 
 \stackrel{\mu}{\ra} & H^1(S - D, \Theta_S( - \log D)).  & & &
\end{array}
$$
Since $(S, Y)$ is  of non-fibered type, from [(2), Proposition  2.1 \cite{STT}], we see that  
$$
H^0(S-D, \Theta_S (- \log D)) =  H^0(S -D, \Theta_S) = \{ 0 \}.
$$ 
Hence we have the following exact sequence:
\begin{equation}\label{eq:local}
\begin{array}{cccccc}
0 \ra 
&  H^1_D(\Theta_S( - \log D)) 
&  \ra 
&  H^1(S, \Theta_S( - \log D))
 & \stackrel{\mu}{\ra} 
 & H^1(S - D, \Theta_S( - \log D)) .\\
\end{array}
\end{equation}

The following theorem is the main statement in this paper.

\begin{Theorem}\label{thm:time}
 Let $(S, Y)$ be  a generalized rational Okamoto-Painlev\'e pair $(S, Y)$ with the condition above. Moreover assume that $D = Y_{red}$ is of additive type. Then  we have 
$$
\dim H^0 (D, \Theta_S(- \log D) \otimes N_{D} ) = 1.
$$
Here  we put $N_D = {\cal O}_S(D)/{\cal O}_S$.  

Since we have 
a natural inclusion 
$$
H^0 (D, \Theta_S(- \log D) \otimes N_{D} ) \hookrightarrow H^1_D(\Theta_S(-\log D) ),
$$ 
we obtain 
$$
\dim H^1_D(\Theta_S(-\log D) ) \geq  1.
$$
\end{Theorem}

This theorem plays an important role to understand the 
Painlev\'e equation related to $(S, Y)$.(cf. \cite{STT}).  
Though we will not investigate the further 
structure of local cohomology here. Instead, we propose 
the following
\begin{Conjecture}\label{conj:local}
Under the same notation and assumption as in Theorem \ref{thm:time},  
$$
H^1_D(\Theta_S(-\log D) ) \simeq H^0(D, \Theta_S(- \log D) \otimes N_{D} ) \simeq \C.
$$
\end{Conjecture}\vspace{2mm}

From the exact sequence (\ref{eq:local}), we see that 
the subspace $ H^1_D(S, \Theta_S(-\log D))$ of 
$H^1(S, \Theta_S(-\log D))$ coincides with the kernel of $\mu$. This implies that:
$$
H^1_D(S, \Theta_S(- \log D)) \simeq 
\left\{ 
\begin{array}{c}
	\mbox{Infinitesimal deformations of $(S, D)$ whose restriction} \\
	\mbox{to $S-D$ induces the trivial deformation}
\end{array}
\right\}.
$$

In [\S 6 \cite{STT}],  
we show  that  any direction 
corresponding
 to a non-zero element of the local cohomology 
group $H^1_D(S, \Theta_S(- \log D))$ induces a 
differential equation (at least locally) by using 
\v{C}ech coboundaries.  

At this moment, we can not  prove Conjecture \ref{conj:local} with  full generality. 
However,  we see that 
the one dimensional 
vector subspace $ H^1(D, \Theta_S(-\log D)\otimes N_D ) $ of $H^1_D(\Theta_S(-\log D) ) \subset H^1(\Theta_S(-\log D) )$ 
really corresponds to the time variable $t$ in the known Painlev\'e equation.  
It is unlikely that we will have more time 
variables, so this gives an  evidence of Conjecture \ref{conj:local}. 
\begin{Remark} {\rm We will consider $(S,Y)$ of multiplicative type later, where the situation is different. (cf. Proposition \ref{prop:multi}).}
\end{Remark}\vspace{2mm}

Let us make preparations for the proof of Theorem \ref{thm:time}.

Recall that 
$$
H^1_D(S, \Theta_S( - \log D)) = \varinjlim {\rm Ext}^1({\cO_{nD}}, \Theta_S( - \log D)) 
$$
where $\cO_{nD} = \cO_S/\cO_S(-n D)$ (cf. [Theorem 2.8, [Gr]]). 

On the other hand, since $\Theta_S( - \log D)$ is a locally free $\cO_S$-module, we see that
$$
{\cal H}om({\cO_{nD}}, \Theta_S( - \log D)) = 0, 
$$
and 
$$
{\cal E}xt^1({\cO_{nD}}, \Theta_S( - \log D)) =  \Theta_S( - \log D) \otimes N_{nD}, 
$$
where $N_{nD} = \cO_S(nD)/\cO_S$.  
By an argument using a spectral sequence, we see that 
$$
H^1_D(S, \Theta_S( - \log D)) = \varinjlim H^0( \Theta_S( - \log D) \otimes N_{nD}) .
$$
Hence we have a natural inclusion
$$
H^0( \Theta_S( - \log D) \otimes N_{D}) \hookrightarrow H^1_D(S, \Theta_S( - \log D)).  
$$

\begin{Lemma}\label{lem:kuma} Let $(S, Y)$ be a generalized rational Okamoto--Painlev\'e pair as above and set $D= Y_{red}$.  
Then we have the following exact sequences :
\begin{equation}\label{eq:thetanormal}
0 \lra \Theta_D \otimes N_D \lra \Theta_S \otimes N_D \lra 
\nu_*(\oplus_{i=1}^r N_{Y_i/S})\otimes N_D \lra 0,
\end{equation}
\begin{equation}\label{eq:thetalog}
0 \lra \nu_*(\oplus_{i=1}^r N_{Y_i/S}) 
 \lra \Theta_S(- \log D) \otimes N_D \lra  \Theta_D \otimes N_D
 \lra 0. 
\end{equation}
Here $\Theta_D$ denotes the tangent sheaf of $D$ and $\nu : \tilde{D}=\coprod^9_{i=1} Y_i \to D$ the normalization map. 
\end{Lemma}

{\it Proof.}  The first exact sequence 
(\ref{eq:thetanormal}) follows from [(1.9), \cite{B-W}].  

Let us consider the following  diagram:  
$$
\begin{array}{ccccccc}
&0&&0&& {\cal \ker \lambda}& \\
 & \downarrow & & \downarrow  & & \downarrow & \\
 &&&&&& \\
0 \lra & 
\Theta_S(- \log D) & 
\lra & 
\Theta_S(- \log D) \otimes \cO_S(D) & 
\lra & 
\Theta_S(- \log D) \otimes N_D & 
\lra 0 \\
 & \downarrow & & \downarrow  & & \quad \downarrow \lambda & \\
0 \lra & 
\Theta_S& 
\lra & 
\Theta_S \otimes \cO_S(D) & 
\lra & 
\Theta_S \otimes N_D & 
\lra 0 \\
 & \downarrow & & \downarrow  & & \downarrow & \\
 & \nu_\ast(\oplus_{i=1}^r N_{Y_i/S}) & \stackrel{\mu}{\lra} &\nu_\ast(\oplus_{i=1}^r N_{Y_i/S}) \otimes N_D  & \lra & \coker \lambda  
  &  \lra 0 \\
   & \downarrow & & \downarrow  & & \downarrow  & \\
&0&&0&& 0 &. \\
\end{array}
$$
By the snake lemma, we obtain the exact sequence 
$$
0 \lra \ker \lambda \lra \nu_\ast(\oplus_{i=1}^r N_{Y_i/S}) \stackrel{\mu}{\lra} \nu_\ast(\oplus_{i=1}^r N_{Y_i/S}) \otimes N_D   \lra  \coker \lambda  
    \lra 0.     
$$
From a local consideration, we see that the map $\mu$ is a zero map, hence
$$
\ker \lambda \simeq \nu_\ast(\oplus_{i=1}^r N_{Y_i/S}),\ \ \mbox{and} \ \ \ \nu_\ast(\oplus_{i=1}^r N_{Y_i/S}) \otimes N_D   \simeq  \coker \lambda.
$$

Moreover since $\im \lambda \simeq \ker [ \Theta_S \otimes N_D \lra \nu_\ast(\oplus_{i=1}^r N_{Y_i/S}) \otimes N_D ]$, from the exact sequence (\ref{eq:thetanormal}), we obtain the exact sequence  (\ref{eq:thetalog}).\qed
\ \\

Note that since $N_{Y_i/S} = \cO_{Y_i}(-2)$, we have 
$$
H^0(\oplus_{i=1}^r N_{Y_i/S}) = \{0 \}, \quad  H^1(\oplus_{i=1}^r N_{Y_i/S} ) \simeq \C^r.
$$
Moreover  one can easily see that
$$
\Theta_D \simeq \nu_\ast(\oplus^r_{i=1} \Theta_{Y_i} ( -t_i)) \simeq 
\nu_\ast((\oplus^r_{i=1} \cO_{Y_i} (2 - t_i) )
$$
where $t_i$ is the number of intersections of $Y_i$ with the other 
$Y_j$.  On the other hand, since $D \cdot Y_i = t_i - 2$ and $\nu$ is a finite morphism, we see that
$$
\begin{array}{cl}
H^0(D,\Theta_D \otimes N_D) & \simeq H^0(D,\nu_\ast (\oplus_{i=1}^r \Theta_{Y_i}(-t_i)) \otimes N_D) \\
 & \simeq H^0(\tilde{D},(\oplus_{i=1}^r \Theta_{Y_i}(-t_i)) \otimes \nu^\ast (N_D)) \\
 & \simeq \oplus_{i=1}^r H^0(Y_i,\Theta_{Y_i}(-t_i) \otimes N_D) \\
 & \simeq \oplus_{i=1}^r H^0(Y_i,\cO_{Y_i}) \\
 & \simeq \C^r.
\end{array}$$

{\it Proof of Theorem \ref{thm:time}}. 

From the exact sequence (\ref{eq:thetalog}), one can obtain 
$$
H^0(\oplus_{i=1}^r N_{Y_i/S}) \lra  
 H^0(\Theta_S(- \log D) \otimes N_D)  \lra H^0( \Theta_D \otimes N_D) \stackrel{\delta}{\lra} H^1(\oplus_{i=1}^r N_{Y_i/S})
$$
where $\delta$ denotes the connected homomorphism. 

From this sequence, the connecting  homomorphism $\delta$  
$$
\delta : H^0(\Theta_D \otimes N_D)  \lra   \oplus_{i=1}^r H^1(N_{Y_i/S}) 
$$
can be identified with  a linear map  $\delta:\C^r \lra \C^r$ and 
we have an isomorphism 
$$
 H^0(D, \Theta_S(- \log D) \otimes N_D ) \simeq \ker \delta. 
$$

Now we state the following proposition, the proof of which is given in \S \ref{sec:proof}.

\begin{Proposition}\label{prop:delta}
Let $(S, Y)$ be a generalized rational Okamoto--Painlev\'{e} pair of additive type and set $D=Y_{red}$. Assume that $D$ is a normal crossing divisor, then we can choose a basis of $H^0(\Theta_D \otimes N_D)$ and $\oplus_{i=1}^r H^1( N_{Y_i} )$ so that the linear map  $ \delta : H^0(\Theta_D \otimes N_D) \ra \oplus_{i=1}^r H^1( N_{Y_i} )$ is represented by the intersection matrix $((Y_i \cdot Y_j))_{1\leq i, j \leq r}$ with respect to these basis.
\end{Proposition}

In this case, it is well-known that the rank of the intersection matrix $((Y_i \cdot Y_j))_{1\leq i, j \leq r}$  is $r-1$. Hence from the Proposition \ref{prop:delta}, we have
$$
\dim_\C H^0(\Theta_S( - \log D) \otimes N_D) = \dim_\C \ker \delta = 1.
$$
\qed

\section{Proof of Proposition \ref{prop:delta}}\label{sec:proof}

In this section, we shall prove Proposition \ref{prop:delta}.

Here we give a detailed proof only for the case of $\tilde{E_7}$. The proof of other cases are similar.

Let $(S,Y)$ be a generalized rational Okamoto--Painlev\'e pair of type $\tilde{E_7}$. Then according to the results in Appendix B of \cite{Sakai}, $(S,Y)$ can be obtained by blowings up 9-points of $\BP^2$ as follows.
\ \\

Let $[x:y:z]$ be the homogeneous coordinates of $\BP^2$.\\
\setlength{\unitlength}{1mm}
\begin{picture}(150,40)(0,0)
	
	\put(11.5,2.5){\line(3,5){15}}
	\put(38.5,2.5){\line(-3,5){15}}
	
	\put(25,25){\circle{2}}
	\put(13,5){\circle*{1}}
	\put(19,15){\circle*{1}}
	
	%%%
	
	\put(29,27){$x=0$}
	\put(13,27){$z=0$}
	
	\put(10,7){$p_4$}
	\put(16,17){$p_5$}
	\put(31,22){$p_1 \gets p_2 \gets p_3 \gets p_6 \gets$}
	\put(41,18){$\gets p_7 \gets p_8 \gets p_9$}

	\put(78,16){$\longleftarrow$}
	
	%%%%%%
	
	\put(95,3){\line(0,1){14}}
	\put(105,13){\line(0,1){14}}
	\put(115,23){\thicklines\line(0,1){14}}
	\put(125,13){\line(0,1){14}}
	\put(135,3){\thicklines\line(0,1){14}}
	\put(93,15){\line(1,0){14}}
	\put(103,25){\line(1,0){24}}
	\put(123,15){\line(1,0){14}}

	\put(95,8){\circle*{1}}
	\put(90,7){$p_9$}
	\put(135,5){\circle*{1}}
	\put(130,4){$p_4$}
	\put(135,10){\circle*{1}}
	\put(130,9){$p_5$}
	
\end{picture}

$$
\begin{array}{l}
	\displaystyle{p_1:(0:1:0)\la p_2:\left(\frac{x}{y} ,\frac{z}{x} \right)=(0,0)\la  p_3:\left(\frac{x}{y} ,\frac{yz}{x^2} \right)=(0,0)\la  p_6:\left(\frac{x}{y} ,\frac{y^2z}{x^3} \right)=(0,1)\la }\\
 	\displaystyle{p_7:\left(\frac{x}{y} ,\frac{y(y^2z-x^3)}{x^4} \right)=(0,0)\la  p_8:\left(\frac{x}{y} ,\frac{y^2(y^2z-x^3)}{x^5} \right)=(0,-s)\la }\\
	\displaystyle{p_9:\left(\frac{x}{y} ,\frac{y(y^2(y^2z-x^3)+sx^5)}{x^6} \right)=(0,-\alpha_0),}\\
	\ \\
	p_4:(0:0:1),\\
	p_5:(0:\alpha_1:1).
\end{array}
$$

Note that there exist three deformation parameters $(\alpha_0,\alpha_1,s)$ of the blowings-up.
Moreover there exists a $\C^\times$-action on the family of surfaces by
$$(\mu ,(\alpha_0,\alpha_1;s,[x:y:z])) \mapsto (\mu^3 \alpha_0,\mu^3 \alpha_1,\mu^2 s;[x:\mu y:\mu^{-2} z]).$$

If we set $\lambda=\alpha_0+\alpha_1$ and $\mu=s/\lambda$, then 
$$
t=\frac{s^3}{\lambda^2},\ a_0=\frac{s^3\alpha_0}{\lambda^3},\ a_1=\frac{s^3\alpha_1}{\lambda^3}, X=x,\ Y=\frac{sy}{\lambda},\ Z=\frac{\lambda^2 z}{s^2}
$$
are invariant under the $\C^\times$-action, and we have the relation
$$
a_0+a_1=t.
$$

Now we can introduce the affine open covering of $S$ and affine coordinates by the explicit blowings-up of $\BP^2$.
The following diagram shows how one can perform blowing-ups and introduce the new coordinates $(x_i,y_i)\ (1 \leq i \leq 13),\ (u_j,v_j)\ (j=14,15,16)$.\\

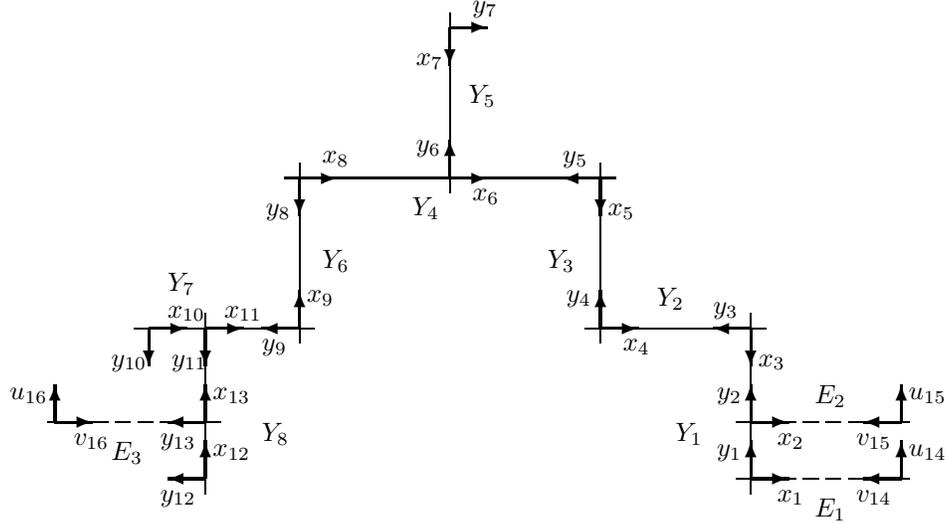
\setlength{\unitlength}{0.5mm}
\begin{picture}(280,140)(-20,0)
%%%%%     E_7     %%%%%
\put(65,10){\thicklines\vector(0,1){10}}
\put(65,10){\thicklines\vector(-1,0){10}}

\put(65, 6){\line(0,1){48}}
\put(80, 20){$Y_8$}

\put(65,50){\thicklines\vector(0,-1){10}}
\put(65,50){\thicklines\vector(1,0){10}}

\put(50,50){\thicklines\vector(0,-1){10}}
\put(50,50){\thicklines\vector(1,0){10}}

\put(46,50){\line(1,0){48}}
\put(55,60){$Y_7$}

\put(90,50){\thicklines\vector(-1,0){10}}
\put(90,50){\thicklines\vector(0,1){10}}

\put(90, 46){\line(0,1){48}}
\put(96,66){$Y_6$}

\put(90,90){\thicklines\vector(0,-1){10}}
\put(90,90){\thicklines\vector(1,0){10}}

\put(86, 90){\line(1,0){88}}
\put(120, 80){$Y_4$}

\put(130,90){\thicklines\vector(0,1){10}}
\put(130,90){\thicklines\vector(1,0){10}}

\put(130,130){\thicklines\vector(0,-1){10}}
\put(130,130){\thicklines\vector(1,0){10}}

\put(130, 86){\line(0,1){48}}
\put(135, 110){$Y_5$} 

\put(170,90){\thicklines\vector(-1,0){10}}
\put(170,90){\thicklines\vector(0,-1){10}}

\put(170, 94){\line(0, -1){48}}
\put(156, 66){$Y_3$}

\put(170,50){\thicklines\vector(0,1){10}}
\put(170,50){\thicklines\vector(1,0){10}}

\put(166, 50){\line(1,0){48}}
\put(185, 56){$Y_2$}

\put(210,50){\thicklines\vector(-1,0){10}}
\put(210,50){\thicklines\vector(0,-1){10}}

\put(210, 54){\line(0, -1){48}}
\put(190, 20){$Y_1$}

%%%

\multiput(206,10)(6,0){8}{\line(1,0){4}}
\put(210,10){\thicklines\vector(1,0){10}}
\put(210,10){\thicklines\vector(0,1){10}}
\put(250,10){\thicklines\vector(-1,0){10}}
\put(250,10){\thicklines\vector(0,1){10}}
\put(227,0){$E_1$}

\multiput(206,25)(6,0){8}{\line(1,0){4}}
\put(210,25){\thicklines\vector(1,0){10}}
\put(210,25){\thicklines\vector(0,1){10}}
\put(250,25){\thicklines\vector(-1,0){10}}
\put(250,25){\thicklines\vector(0,1){10}}
\put(227,30){$E_2$}

\multiput(69,25)(-6,0){8}{\line(-1,0){4}}
\put(65,25){\thicklines\vector(-1,0){10}}
\put(65,25){\thicklines\vector(0,1){10}}
\put(25,25){\thicklines\vector(1,0){10}}
\put(25,25){\thicklines\vector(0,1){10}}
\put(40,15){$E_3$}
%%%%%%%%%

\put(217,4){$x_1$}
\put(201,16){$y_1$}
\put(217,19){$x_2$}
\put(201,31){$y_2$}
\put(212,40){$x_3$}
\put(200,54){$y_3$}
\put(176,44){$x_4$}
\put(161,57){$y_4$}
\put(172,80){$x_5$}
\put(160,94){$y_5$}
\put(136,84){$x_6$}
\put(121,97){$y_6$}
\put(121,120){$x_7$}
\put(136,134){$y_7$}
\put(96,94){$x_8$}
\put(81,80){$y_8$}
\put(92,57){$x_9$}
\put(80,44){$y_9$}
\put(55,53){$x_{10}$}
\put(40,40){$y_{10}$}
\put(70,53){$x_{11}$}
\put(56,40){$y_{11}$}
\put(67,16){$x_{12}$}
\put(53,4){$y_{12}$}
\put(67,31){$x_{13}$}
\put(53,19){$y_{13}$}
\put(252,16){$u_{14}$}
\put(238,4){$v_{14}$}
\put(252,31){$u_{15}$}
\put(238,19){$v_{15}$}
\put(13,31){$u_{16}$}
\put(30,19){$v_{16}$}

\end{picture}
$$
\begin{array}{c}
	E_i \simeq \BP^1,\ E_i^2=-1,\\
	Y_i \simeq \BP^1,\ Y_i^2=-2.
\end{array}
$$

\begin{eqnarray}
	U_i & = & \Spec \C \left[x_i,y_i\right] \ \cong \ {\C}^2\ \ (i=1,2,\cdots,7).\nonumber\\
	U_8 & = & \Spec \C \left[x_8,y_8,\frac{1}{1+x_8}\right] \ \cong \ {\C}^2-\{1+x_8=0\}.\nonumber\\
	U_9 & = & \Spec \C \left[x_9,y_9,\frac{1}{1 + {x_9}^2{y_9}}\right] \ \cong \ {\C}^2-\{1 + {x_9}^2{y_9}=0\}.\nonumber\\
	U_{10} & = & \Spec \C \left[x_{10},y_{10},\frac{1}{1 +x_{10}{y_{10}}^2}\right] \ \cong \ {\C}^2-\{1 +x_{10}{y_{10}}^2=0\}.\nonumber\\
	U_{11} & = & \Spec \C \left[x_{11},y_{11},\frac{1}{1 - t{x_{11}}^2{y_{11}}^2 + {x_{11}}^3{y_{11}}^2}\right] \ \cong \ {\C}^2-\{1 - t{x_{11}}^2{y_{11}}^2 + {x_{11}}^3{y_{11}}^2=0\}.\nonumber\\
	U_{12} & = & \Spec \C \left[x_{12},y_{12},\frac{1}{1 - t{y_{12}}^2 + x_{12}{y_{12}}^3}\right] \ \cong \ {\C}^2-\{1 - t{y_{12}}^2 + x_{12}{y_{12}}^3=0\}.\nonumber\\
	U_{13} & = & \Spec \C \left[x_{13},y_{13},\frac{1}{-1 + t{x_{13}}^2{y_{13}}^2 + a_0{x_{13}}^3{y_{13}}^3 - {x_{13}}^4{y_{13}}^3}\right]\nonumber\\
	   & \cong & {\C}^2-\{-1 + t{x_{13}}^2{y_{13}}^2 + a_0{x_{13}}^3{y_{13}}^3 - {x_{13}}^4{y_{13}}^3=0\}.\nonumber\\
	U_j & = & \Spec \C \left[u_j,v_j\right] \ \cong \ {\C}^2\ \ (j=14,15).\nonumber\\
	U_{16}& = & \Spec \C \left[u_{16},v_{16},\frac{1}{-1+tu_{16}^2 + a_0u_{16}^3 - u_{16}^4v_{16}}\right] \ \cong \ {\C}^2-\{-1+tu_{16}^2 + a_0u_{16}^3 - u_{16}^4v_3=0\}.\nonumber
\end{eqnarray}
$$
\begin{array}{ll}
	Y_1=\{x_1=0,x_2=0,y_3=0 \},& Y_2=\{x_3=0,y_4=0 \},\vspace{2mm}\\
	Y_3=\{x_4=0,y_5=0 \},& Y_4=\{x_5=0,y_6=0,y_8=0 \},\vspace{2mm}\\
	Y_5=\{x_6=0,y_7=0 \},& Y_6=\{x_8=0,y_9=0 \},\vspace{2mm}\\
	Y_7=\{x_9=0,y_{10}=0,y_{11}=0 \},& Y_8=\{x_{11}=0,y_{12}=0,y_{13}=0 \}.
\end{array}
$$
$$
\begin{array}{c}
\displaystyle{S=\bigcup_{i=1}^{16}U_i.}\vspace{2mm}\\
\displaystyle{Y=Y_1+2Y_2+3Y_3+4Y_4+2Y_5+3Y_6+2Y_7+Y_8,\ \ \ D=\sum^8_{i=1}Y_i.}\vspace{2mm}\\
S-Y=U_{14} \cup U_{15} \cup U_{16}.
\end{array}
$$

More explicitly, new coordinates can be given by the following formula
$$
\begin{array}{c}
\left\{
\begin{array}{l}
\displaystyle{x_1=\frac{X}{Y}}\vspace{1mm}\\
\displaystyle{y_1=\frac{Y}{Z}}
\end{array}
\right.
\hspace{2cm}
\left\{
\begin{array}{l}
\displaystyle{x_2=\frac{X}{Y-a_1 Z}}\vspace{1mm}\\
\displaystyle{y_2=\frac{Y-a_1 Z}{Z}}
\end{array}
\right.
\hspace{2cm}
\left\{
\begin{array}{l}
\displaystyle{x_3 =\frac{Z}{Y}}\vspace{1mm}\\
\displaystyle{y_3 =\frac{X}{Z}}
\end{array}
\right.
\vspace{2mm}\\
\left\{
\begin{array}{l}
\displaystyle{x_4 =\frac{Z}{X}}\vspace{1mm}\\
\displaystyle{y_4 =\frac{X^2}{YZ}}
\end{array}
\right.
\hspace{2cm}
\left\{
\begin{array}{l}
\displaystyle{x_5 =\frac{YZ}{X^2}}\vspace{1mm}\\
\vspace{3mm}
\displaystyle{y_5 =\frac{X^3}{Y^2Z}}
\end{array}
\right.
\hspace{2cm}
\left\{
\begin{array}{l}
\displaystyle{x_6 =\frac{Y^2Z}{X^3}}\vspace{1mm}\\
\displaystyle{y_6 =\frac{X}{Y}}
\end{array}
\right.
\vspace{2mm}\\
\left\{
\begin{array}{l}
\displaystyle{x_7 =\frac{Y}{X}}\vspace{1mm}\\
\displaystyle{y_7 =\frac{Z}{X}}
\end{array}
\right.
\hspace{1cm}
\left\{
\begin{array}{l}
\displaystyle{x_8 =\frac{Y^2Z - X^3}{X^3}}\vspace{1mm}\\
\displaystyle{y_8 =\frac{X^4}{Y(Y^2Z-X^3)}}
\end{array}
\right.
\hspace{1cm}
\left\{
\begin{array}{l}
\displaystyle{x_9 =\frac{Y(Y^2Z-X^3)}{X^4}}\vspace{1mm}\\
\displaystyle{y_9 =\frac{X^5}{Y^2(Y^2Z-X^3)}}
\end{array}
\right.
\vspace{2mm}\\
\left\{
\begin{array}{l}
\displaystyle{x_{10} =\frac{Y^2(Y^2Z-X^3)}{X^5}}\vspace{1mm}\\
\displaystyle{y_{10} =\frac{X}{Y}}
\end{array}
\right.
\hspace{2cm}
\left\{
\begin{array}{l}
\displaystyle{x_{11} =\frac{Y^2(Y^2Z-X^3)+tX^5}{X^5}}\vspace{1mm}\\
\displaystyle{y_{11} =\frac{X^6}{Y(Y^2(Y^2Z-X^3)+tX^5)}}
\end{array}
\right.
\vspace{2mm}\\
\left\{
\begin{array}{l}
\displaystyle{x_{12} =\frac{Y(Y^2(Y^2Z-X^3)+tX^5)}{X^6}}\vspace{1mm}\\
\displaystyle{y_{12} =\frac{X}{Y}}
\end{array}
\right.
\left\{
\begin{array}{l}
\displaystyle{x_{13} =\frac{Y(Y^2(Y^2Z-X^3)+tX^5)+a_0 X^6}{X^6}}\vspace{1mm}\\
\displaystyle{y_{13} =\frac{X^7}{Y(Y(Y^2(Y^2Z-X^3)+tX^5)+a_0 X^6)}}
\end{array}
\right.
\vspace{2mm}\\
\left\{
\begin{array}{l}
\displaystyle{u_{14} =\frac{X}{Z}}\vspace{1mm}\\
\displaystyle{v_{14} =\frac{Y}{X}}
\end{array}
\right.
\ \ 
\left\{
\begin{array}{l}
\displaystyle{u_{15} =\frac{X}{Z}}\vspace{1mm}\\
\displaystyle{v_{15} =\frac{Y-a_1 Z}{X}}
\end{array}
\right.
\ \ 
\left\{
\begin{array}{l}
\displaystyle{u_{16} =\frac{X}{Y}}\vspace{1mm}\\
\displaystyle{v_{16} =\frac{Y(Y(Y^2(Y^2Z-X^3)+tX^5)+a_0 X^6)}{X^7}}.
\end{array}
\right.
\end{array}
$$

\ \\
From these formula, we can determine the coordinate transformation between $(x_i,y_i)$'s and $(u_j,v_j)$'s.\\
For later use, we need only the coordinate transformations near each component $Y_i$.
Here we will list up the coordinate transformations only for a neighborhood of each $Y_i$.
$$
\hspace*{-2mm}
\begin{array}{ll}
Y_1:
\left\{
\begin{array}{l}
	\displaystyle{x_1=\frac{x_2y_2}{a_1 + y_2}}\\
	\displaystyle{y_1=a_1 + y_2}
\end{array}
\right.
\left\{
\begin{array}{l}
	\displaystyle{x_1=x_3y_3}\\
	\displaystyle{y_1=\frac{1}{x_3}}
\end{array}
\right.
&
Y_5:\left\{
\begin{array}{l}
	\displaystyle{x_6={x_7}^2y_7}\\
	\displaystyle{y_6=\frac{1}{x_7}}
\end{array}
\right.
\vspace{2mm}\\
Y_2:\left\{
\begin{array}{l}
	\displaystyle{x_3={x_4}^2y_4}\\
	\displaystyle{y_3=\frac{1}{x_4}}
\end{array}
\right.
&
Y_6:\left\{
\begin{array}{l}
	\displaystyle{x_8={x_9}^2y_9}\\
	\displaystyle{y_8=\frac{1}{x_9}}
\end{array}
\right.
\vspace{2mm}\\
Y_3:\left\{
\begin{array}{l}
	\displaystyle{x_4={x_5}^2y_5}\\
	\displaystyle{y_4=\frac{1}{x_5}}
\end{array}
\right.
&
Y_7:\left\{
\begin{array}{l}
	\displaystyle{x_9=x_{10}y_{10}}\\
	\displaystyle{y_9=\frac{1}{x_{10}}}
\end{array}
\right.
\left\{
\begin{array}{l}
	\displaystyle{x_9=x_{11} (-t+x_{11}) y_{11}}\\
	\displaystyle{y_9=\frac{1}{-t+x_{11}}}
\end{array}
\right.
\vspace{2mm}\\
Y_4:\left\{
\begin{array}{l}
	\displaystyle{x_5=x_6y_6}\\
	\displaystyle{y_5=\frac{1}{x_6}}
\end{array}
\right.
\left\{
\begin{array}{l}
	\displaystyle{x_5=x_8 (1+x_8) y_8}\\
	\displaystyle{y_5=\frac{1}{1+x_8}}
\end{array}
\right.
&
Y_8:\left\{
\begin{array}{l}
	\displaystyle{x_{11}=x_{12}y_{12}}\\
	\displaystyle{y_{11}=\frac{1}{x_{12}}}
\end{array}
\right.
\left\{
\begin{array}{l}
	\displaystyle{x_{11}=x_{13} (-a_0+x_{13}) y_{13}}\\
	\displaystyle{y_{11}=\frac{1}{-a_0+x_{13}}}.
\end{array}
\right.
\end{array}
$$

Let us consider the sheaf $\Theta_S(- \log D)$ 
and the sheaf exact sequence 
$$
0 \lra \nu_\ast(\oplus_{i=1}^8 N_{Y_i/S}) \lra 
\Theta_S(- \log D) \otimes N_D \lra \Theta_D \otimes N_D \lra 0.
$$
We will analyse the edge homomorphism
\begin{equation}\label{eq:intes}
\delta : H^0(\Theta_D \otimes N_D) \lra H^1(\nu_\ast(\oplus_{i=1}^8 N_{Y_i /S}))
\end{equation}
by using the \v{C}ech cocycles.

Noting that $\nu$ is a finite morphism, and
$ \Theta_D \otimes N_D 
\simeq \nu_\ast (\oplus_{i=1}^8 \Theta_{Y_i}(-t_i)) \otimes N_D$  where $t_i$ is the number of intersections of $Y_i$ with other components, we see that
\begin{equation}\label{eq:kamo}
H^0(D,\Theta_D \otimes N_D) \simeq H^0(D,\nu_\ast (\oplus_{i=1}^8 \Theta_{Y_i}(-t_i)) \otimes N_D) \simeq \C^8.
\end{equation}
For each $i\ (1 \le i \le 8)$, we introduce a generator $\theta_i$ of the cohomology group in (\ref{eq:kamo}) as follows.\\

\begin{center}
\begin{tabular}{|c|l|}\hline
&\vspace{-2mm}\\
$\theta_1$ & 
$\left\{
\begin{array}{l}
\theta_1^1=\frac{-a_1+y_1}{x_1}\frac{\partial}{\partial y_1}\ \mbox{on} \ U_1\cap Y_1 ,\ \ 
\theta_1^2=\frac{a_1+y_2}{x_2}\frac{\partial}{\partial y_2}\ \mbox{on} \ U_2\cap Y_1 ,\vspace{1mm}\\
\theta_1^3=\frac{-1+a_1x_3}{y_3}\frac{\partial}{\partial x_3}\ \mbox{on} \ U_3\cap Y_1  
\end{array}
\right\}$\\
&\vspace{-2mm}\\
\hline
&\vspace{-2mm}\\
$\theta_2 $& 
$\left\{\theta_2^3=\frac{1}{x_3}\frac{\partial}{\partial y_3}\ \mbox{on} \ U_3\cap Y_2 ,\ \ 
\theta_2^4=-\frac{1}{y_4}\frac{\partial}{\partial x_4}\ \mbox{on} \ U_4\cap Y_2  \right\}$\\
&\vspace{-2mm}\\
\hline
&\vspace{-2mm}\\
$\theta_3$ & 
$\left\{\theta_3^4=\frac{1}{x_4}\frac{\partial}{\partial y_4}\ \mbox{on} \ U_4\cap Y_3 ,\ \ 
\theta_3^5=-\frac{1}{y_5}\frac{\partial}{\partial x_5}\ \mbox{on} \ U_5\cap Y_3  \right\}$\\
&\vspace{-2mm}\\
\hline
&\vspace{-2mm}\\
$\theta_4$ & 
$\left\{
\begin{array}{l}
\theta_4^5=\frac{1-y_5}{x_5}\frac{\partial}{\partial y_5}\ \mbox{on} \ U_5\cap Y_4 ,\ \ 
\theta_4^6=\frac{1-x_6}{y_6}\frac{\partial}{\partial x_6}\ \mbox{on} \ U_6\cap Y_4 ,\vspace{1mm}\\
\theta_4^8=-\frac{1}{y_8}\frac{\partial}{\partial x_8}\ \mbox{on} \ U_8\cap Y_4 
\end{array} 
\right\}$\\
&\vspace{-2mm}\\
\hline
&\vspace{-2mm}\\
$\theta_5$ & 
$\left\{\theta_5^6=-\frac{1}{x_6}\frac{\partial}{\partial y_6}\ \mbox{on} \ U_6\cap Y_5 ,\ \ 
\theta_5^7=\frac{1}{y_7}\frac{\partial}{\partial x_7}\ \mbox{on} \ U_7\cap Y_5  \right\}$\\
&\vspace{-2mm}\\
\hline
&\vspace{-2mm}\\
$\theta_6$ & 
$\left\{\theta_6^8=\frac{1}{x_8}\frac{\partial}{\partial y_8}\ \mbox{on} \ U_8\cap Y_6 ,\ \ 
\theta_6^9=-\frac{1}{y_9}\frac{\partial}{\partial x_9}\ \mbox{on} \ U_8\cap Y_6  \right\}$\\
&\vspace{-2mm}\\
\hline
&\vspace{-2mm}\\
$\theta_7$ & 
$\left\{
\begin{array}{l}
\theta_7^9=\frac{1+ty_9}{x_9}\frac{\partial}{\partial y_9}\ \mbox{on} \ U_9\cap Y_7 ,\ \ 
\theta_7^{10}=-\frac{t+x_{10}}{y_{10}}\frac{\partial}{\partial x_{10}}\ \mbox{on} \ U_{10}\cap Y_7 ,\vspace{1mm}\\
\theta_7^{11}=-\frac{1}{y_{11}}\frac{\partial}{\partial x_{11}}\ \mbox{on} \ U_{11}\cap Y_7
\end{array}
\right\}$\\
&\vspace{-2mm}\\
\hline
&\vspace{-2mm}\\
$\theta_8$ & 
$\left\{
\begin{array}{l}
\theta_8^{11}=\frac{1+a_0y_{11}}{x_{11}}\frac{\partial}{\partial y_{11}}\ \mbox{on} \ U_{11}\cap Y_8  ,\ \ 
\theta_8^{12}=-\frac{a_0+x_{12}}{y_{12}}\frac{\partial}{\partial x_{12}}\ \mbox{on} \ U_{12}\cap Y_8,\vspace{1mm}\\
\theta_8^{13}=-\frac{1}{y_{13}}\frac{\partial}{\partial x_{13}}\ \mbox{on} \ U_{13}\cap Y_8
\end{array}
\right\}$\vspace{-2mm}\\
&\\
\hline
\end{tabular}
\end{center}

On the other hand, for each $i\ (1 \le i \le 8)$, we have a generator $\eta_i \in H^1(Y_i,N_{Y_i/S})$ as follows.
\begin{center}
\begin{tabular}{|c|l|}\hline
& \vspace{-2mm}\\
$\eta_1$ & 
$\left\{\eta_1^{12}=0\ \mbox{on} \ U_1\cap \ U_2\cap Y_1 ,\ \ \eta_1^{13}=\frac{1}{x_3}\frac{\partial}{\partial y_3}\ \mbox{on} \ U_1\cap \ U_3\cap Y_1  \right\}$\\
&\vspace{-2mm}\\
\hline
&\vspace{-2mm}\\
$\eta_2$ & 
$\left\{\eta_2^{34}=\frac{1}{x_4}\frac{\partial}{\partial y_4}\ \mbox{on} \ U_3\cap \ U_4\cap Y_2 \right\}$\\
&\vspace{-2mm}\\
\hline
&\vspace{-2mm}\\
$\eta_3$ & 
$\left\{\eta_3^{45}=\frac{1}{x_5}\frac{\partial}{\partial y_5}\ \mbox{on} \ U_4\cap \ U_5\cap Y_3 \right\}$\\
&\vspace{-2mm}\\
\hline
&\vspace{-2mm}\\
$\eta_4$ & 
$\left\{\eta_4^{56}=0\ \mbox{on} \ U_5\cap \ U_6\cap Y_4 ,\ \ \eta_4^{58}=\frac{1}{x_8}\frac{\partial}{\partial y_8}\ \mbox{on} \ U_5\cap \ U_8\cap Y_4 \right\}$\\
&\vspace{-2mm}\\
\hline
&\vspace{-2mm}\\
$\eta_5$ & 
$\left\{\eta_5^{67}=-\frac{1}{x_7}\frac{\partial}{\partial y_7}\ \mbox{on} \ U_6\cap \ U_7\cap Y_5 \right\}$\\
&\vspace{-2mm}\\
\hline
&\vspace{-2mm}\\
$\eta_6$ & 
$\left\{\eta_6^{89}=\frac{1}{x_9}\frac{\partial}{\partial y_9}\ \mbox{on} \ U_8\cap \ U_9\cap Y_6  \right\}$\\
&\vspace{-2mm}\\
\hline
&\vspace{-2mm}\\
$\eta_7$ & 
$\left\{\eta_7^{9\, 10}=0\ \mbox{on} \ U_9\cap \ U_{10}\cap Y_7,\ \ \eta_7^{9\, 11}=\frac{1}{x_{11}}\frac{\partial}{\partial y_{11}}\ \mbox{on} \ U_9\cap \ U_{11}\cap Y_7 \right\}$\\
&\vspace{-2mm}\\
\hline
&\vspace{-2mm}\\
$\eta_8$ & 
$\left\{\eta_8^{11\, 13}=\frac{1}{x_{13}}\frac{\partial}{\partial y_{13}}\ \mbox{on} \ U_{11}\cap \ U_{13}\cap Y_8,\ \ \eta_8^{12\, 13}=0 \ \mbox{on} \ U_{12}\cap \ U_{13}\cap Y_8\right\}$\vspace{-2mm}\\
&\\
\hline
\end{tabular}
\end{center}

\ \\
We take $\{\theta_i\}$ and $\{\eta_i\}$ as basis of $H^0(\Theta_D \otimes N_D)$ and $\oplus_{i=1}^8H^1(N_{Y_i /M})$ respectively.\\
\ \\

By using these bases, we compute the matrix representing the connecting homomorphism $\delta$.

For that purpose, let us lift 0-cocycle $\theta_1$ to 0-cochains of $\Theta_S(- \log D) \otimes N_D$ as 
$$
\tilde{\theta_1^1}=\frac{-a_1+y_1}{x_1}\frac{\partial}{\partial y_1}\ {\rm on}\ U_1,\hspace{5mm}
\tilde{\theta_1^2}=\frac{a_1+y_2}{x_2}\frac{\partial}{\partial y_2}\ {\rm on}\ U_2,\hspace{5mm}
\tilde{\theta_1^3}=\frac{-1+a_1x_3}{y_3}\frac{\partial}{\partial x_3}\ {\rm on}\ U_3,
$$
$$\tilde{\theta_1^i}=0\ {\rm on}\ U_i\ (i=4,5,\cdots,16).$$
Other 0-cocycles can be lifted in a similar way.\\

We first compute $\delta(\theta_1)$.\\

From the definition of $\delta$, we have $\delta(\theta_1)=\{ \delta(\theta_1)_{ij} \ {\rm on}\ U_i \cap U_j \cap D \}$ with
$$
\begin{array}{l}
	\left\{
	\begin{array}{l}
		\displaystyle{\delta(\theta_1)_{12}
		=(\tilde{\theta_1^2}-\tilde{\theta_1^1})|_{Y_1}
		=\left.\left(\frac{a_1+y_2}{x_2}\frac{\partial}{\partial y_2}-\frac{-a_1+y_1}{x_1}\frac{\partial}{\partial y_1}\right)\right|_{Y_1}
		=\frac{a_1}{y_2}\frac{\partial}{\partial x_2} }\\
		
		\displaystyle{\delta(\theta_1)_{13}
		=(\tilde{\theta_1^3}-\tilde{\theta_1^1})|_{Y_1}
		=\left.\left( \frac{-1+a_1x_3}{y_3}\frac{\partial}{\partial x_3}-\frac{-a_1+y_1}{x_1}\frac{\partial}{\partial y_1}\right)\right|_{Y_1}
		=\frac{-1+a_1x_3}{x_3}\frac{\partial}{\partial y_3} }
	\end{array}
	\right.\\
	\displaystyle{\delta(\theta_1)_{34}
	=(\tilde{\theta_1^4}-\tilde{\theta_1^3})|_{Y_2}
	=\left.\left(0-\frac{-1+a_1x_3}{y_3}\frac{\partial}{\partial x_3}\right)\right|_{Y_2}
	=\frac{1}{x_4}\frac{\partial}{\partial y_4}
	=\eta_2^{34}}\\
\end{array}
$$
Other $\delta(\theta_1)_{ij}$'s are zero.\\

Obviously $\{\delta(\theta_1)_{34} \}=\{\eta_2^{34} \}=\eta_2$.

In order to see that $ \{ \delta(\theta_1)_{12},\delta(\theta_1)_{13} \} =-2\eta_1$,
we set $\tau=\{\tau_1=-\frac{\partial}{\partial x_1},\tau_2=-\frac{\partial}{\partial x_2},\tau_3=a_1\frac{\partial}{\partial y_3}\}\in C^0(N_{Y_1 /S})$.
$$
\left\{
\begin{array}{l}
	\displaystyle{\delta(\theta_1)_{12}+2\eta_1^{12}=\frac{a_1}{y_2}\frac{\partial}{\partial x_2}=\tau_2 -\tau_1  }\\
	\displaystyle{\delta(\theta_1)_{13}+2\eta_1^{13}=\frac{1+a_1x_3}{x_3}\frac{\partial}{\partial y_3}=\tau_3 -\tau_1  }
\end{array}
\right.
$$
This implies $ \{ \delta(\theta_1)_{12},\delta(\theta_1)_{13} \}+2\eta_1=\delta\tau$.
Then we have $ \delta(\theta_1)=-2\eta_1+\eta_2$.\\

Other $\delta(\theta_i)$'s can be treated in a similar way. In what follows, we just list up the results of computations.\\

$\circ\  \delta(\theta_2)=\eta_1 -2\eta_2 +\eta_3$
$$
\begin{array}{l}
	\left\{
	\begin{array}{l}
		\displaystyle{\delta(\theta_2)_{12}=\tilde{\theta_2^2}-\tilde{\theta_2^1}}|_{Y_1}=0-0=\eta_1^{12}\\
		\displaystyle{\delta(\theta_2)_{13}=\tilde{\theta_2^3}-\tilde{\theta_2^1}|_{Y_1}=\frac{1}{x_3}\frac{\partial}{\partial y_3}-0=\eta_1^{13}}
	\end{array}
	\right.\\
	\displaystyle{\delta(\theta_2)_{34}=\tilde{\theta_2^4}-\tilde{\theta_2^3}|_{Y_2}=-\frac{1}{y_4}\frac{\partial}{\partial x_4}-\frac{1}{x_3}\frac{\partial}{\partial y_3}=-2\frac{1}{x_4}\frac{\partial}{\partial y_4}=-2\eta_2^{34} }\\
	\displaystyle{\delta(\theta_2)_{45}=\tilde{\theta_2^5}-\tilde{\theta_2^4}|_{Y_3}=0-\left(-\frac{1}{y_4}\frac{\partial}{\partial x_4}\right)=\frac{1}{x_5}\frac{\partial}{\partial y_5}=\eta_3^{45} }
\end{array}
$$
\ \\
$\circ\  \delta(\theta_3)=\eta_2 -2\eta_3 +\eta_4$
$$
\begin{array}{l}
	\displaystyle{\delta(\theta_3)_{34}=\tilde{\theta_3^4}-\tilde{\theta_3^3}|_{Y_2}=\frac{1}{x_4}\frac{\partial}{\partial y_4}-0=\frac{1}{x_4}\frac{\partial}{\partial y_4}=\eta_2^{34} }\\
	\displaystyle{\delta(\theta_3)_{45}=\tilde{\theta_3^5}-\tilde{\theta_3^4}|_{Y_3}=-\frac{1}{y_5}\frac{\partial}{\partial x_5}-\frac{1}{x_4}\frac{\partial}{\partial y_4}=-2\frac{1}{x_5}\frac{\partial}{\partial y_5}=-2\eta_3^{45} }\\
	\left\{
	\begin{array}{l}
		\displaystyle{\delta(\theta_3)_{56}=\tilde{\theta_3^6}-\tilde{\theta_3^5}|_{Y_4}=0-\left(-\frac{1}{y_5}\frac{\partial}{\partial x_5}\right)=\frac{\partial}{\partial y_6} }\\
		\displaystyle{\delta(\theta_3)_{58}=\tilde{\theta_3^8}-\tilde{\theta_3^5}|_{Y_4}=0-\left(-\frac{1}{y_5}\frac{\partial}{\partial x_5}\right)=\frac{1}{x_8}\frac{\partial}{\partial y_8} }
	\end{array}
	\right.
\end{array}
$$
\ \\
\hspace{1cm} Set $\{\tau_5=0,\tau_6=\frac{\partial}{\partial y_6},\tau_8=0 \}\in C^0(N_{Y_4 /S})$.
$$
\left\{
\begin{array}{l}
	\displaystyle{\delta(\theta_3)_{56}-\eta_4^{56}=\frac{\partial}{\partial y_6}-0=\tau_6 -\tau_5  }\\
	\displaystyle{\delta(\theta_3)_{58}-\eta_4^{58}=\frac{1}{x_8}\frac{\partial}{\partial y_8}-\frac{1}{x_8}\frac{\partial}{\partial y_8}=0=\tau_8 -\tau_5 }
\end{array}
\right.
$$
\hspace{1cm} Thus we have $ \{ \delta(\theta_3)_{56},\delta(\theta_3)_{58} \} =\eta_4$.\\
\ \\

$\circ\  \delta(\theta_4)=\eta_3 -2\eta_4 +\eta_5 +\eta_6$
$$
\begin{array}{l}
	\displaystyle{\delta(\theta_4)_{45}=\tilde{\theta_4^5}-\tilde{\theta_4^4}|_{Y_3}=\frac{1-y_5}{x_5}\frac{\partial}{\partial y_5}-0=\frac{1}{x_5}\frac{\partial}{\partial y_5}=\eta_3^{45} }\\
	\left\{
	\begin{array}{l}
		\displaystyle{\delta(\theta_4)_{56}=\tilde{\theta_4^6}-\tilde{\theta_4^5}|_{Y_4}=\frac{1-x_6}{y_6}\frac{\partial}{\partial x_6}-\frac{1-y_5}{x_5}\frac{\partial}{\partial y_5}=\frac{1-x_6}{x_6}\frac{\partial}{\partial y_6} }\\
		\displaystyle{\delta(\theta_4)_{58}=\tilde{\theta_4^8}-\tilde{\theta_4^5}|_{Y_4}=-\frac{1}{y_8}\frac{\partial}{\partial x_8}-\frac{1-y_5}{x_5}\frac{\partial}{\partial y_5}=-\frac{1+2x_8}{x_8(1+x_8)}\frac{\partial}{\partial y_8} }
	\end{array}
	\right.\\
	\displaystyle{\delta(\theta_4)_{67}=\tilde{\theta_4^7}-\tilde{\theta_4^6}|_{Y_5}=0-\frac{1-x_6}{y_6}\frac{\partial}{\partial x_6}=-\frac{1}{x_7}\frac{\partial}{\partial y_7}=\eta_5^{67} }\\
	\displaystyle{\delta(\theta_4)_{89}=\tilde{\theta_4^9}-\tilde{\theta_4^8}|_{Y_6}=0-\left(-\frac{1}{y_8}\frac{\partial}{\partial x_8}\right)=\frac{1}{x_9}\frac{\partial}{\partial y_9}=\eta_6^{89} }
\end{array}
$$
\ \\
\hspace{1cm} Set $\{\tau_5=-\frac{\partial}{\partial x_5},\tau_6=-\frac{\partial}{\partial y_6},\tau_8=0\}\in C^0(N_{Y_4 /S})$.
$$
\left\{
\begin{array}{l}
	\displaystyle{\delta(\theta_4)_{56}+2\eta_4^{56}=\frac{1-x_6}{x_6}\frac{\partial}{\partial y_6}=\tau_6 -\tau_5  }\\
	\displaystyle{\delta(\theta_4)_{58}+2\eta_4^{56}=\frac{1}{x_8(1+x_8)}\frac{\partial}{\partial y_8}=\tau_8 -\tau_5  }
\end{array}
\right.
$$
\hspace{1cm} Thus we have $ \{ \delta(\theta_4)_{56},\delta(\theta_4)_{58} \} =-2\eta_4$.\\
\ \\

$\circ\  \delta(\theta_5)=\eta_4 -2\eta_5$
$$
\begin{array}{l}
	\left\{
	\begin{array}{l}
		\displaystyle{\delta(\theta_5)_{56}=\tilde{\theta_5^6}-\tilde{\theta_5^5}|_{Y_4}=-\frac{1}{x_6}\frac{\partial}{\partial y_6}-0=-\frac{1}{x_6}\frac{\partial}{\partial y_6} }\\
		\displaystyle{\delta(\theta_5)_{58}=\tilde{\theta_5^8}-\tilde{\theta_5^5}|_{Y_4}=0-0=0 }
	\end{array}
	\right.\\
	\displaystyle{\delta(\theta_5)_{67}=\tilde{\theta_5^7}-\tilde{\theta_5^6}|_{Y_5}=\frac{1}{y_7}\frac{\partial}{\partial x_7}-\left(-\frac{1}{x_6}\frac{\partial}{\partial y_6}\right)=2\frac{1}{x_7}\frac{\partial}{\partial y_7}=-2\eta_5^{67} }
\end{array}
$$
\hspace{1cm} Set $\{\tau_5=\frac{\partial}{\partial x_5},\tau_6=0,\tau_8=0\}\in C^0(N_{Y_4 /S})$.
$$
\left\{
\begin{array}{l}
	\displaystyle{\delta(\theta_5)_{56}-\eta_4^{56}=-\frac{1}{x_6}\frac{\partial}{\partial y_6}=\tau_6 -\tau_5 }\\
	\displaystyle{\delta(\theta_5)_{58}-\eta_4^{58}=-\frac{1}{x_8(1+x_8)}\frac{\partial}{\partial y_8}=\tau_8 -\tau_5 }
\end{array}
\right.
$$
\hspace{1cm} Thus we have $ \{ \delta(\theta_5)_{56},\delta(\theta_5)_{58} \} =\eta_4$.\\
\ \\

$\circ\  \delta(\theta_6)=\eta_4 -2\eta_6 +\eta_7$
$$
\begin{array}{l}
	\left\{
	\begin{array}{l}
		\displaystyle{\delta(\theta_6)_{56}=\tilde{\theta_6^6}-\tilde{\theta_6^5}|_{Y_4}=0-0=\eta_4^{56} }\\
		\displaystyle{\delta(\theta_6)_{58}=\tilde{\theta_6^8}-\tilde{\theta_6^5}|_{Y_4}=\frac{1}{x_8}\frac{\partial}{\partial y_8}-0=\eta_4^{58} }
	\end{array}
	\right.\\
	\displaystyle{\delta(\theta_6)_{89}=\tilde{\theta_6^9}-\tilde{\theta_6^8|_{Y_6}}=-\frac{1}{y_9}\frac{\partial}{\partial x_9}-\frac{1}{x_8}\frac{\partial}{\partial y_8}=-2\frac{1}{x_9}\frac{\partial}{\partial y_9}=-2\eta_6^{89} }\\
	\left\{
	\begin{array}{l}
		\displaystyle{\delta(\theta_6)_{9\, 10}=\tilde{\theta_6^{10}}-\tilde{\theta_6^9}|_{Y_7}=0-\left(-\frac{1}{y_9}\frac{\partial}{\partial x_9}\right)=\frac{\partial}{\partial y_{10}} }\\
		\displaystyle{\delta(\theta_6)_{9\, 11}=\tilde{\theta_6^{11}}-\tilde{\theta_6^9}|_{Y_7}=0-\left(-\frac{1}{y_9}\frac{\partial}{\partial x_9}\right)=\frac{1}{x_{11}}\frac{\partial}{\partial y_{11}} }
	\end{array}
	\right.
\end{array}
$$
\ \\
\hspace{1cm} Set $\{\tau_9=0,\tau_{10}=\frac{\partial}{\partial y_{10}},\tau_{11}=0\}\in C^0(N_{Y_7 /S})$.
$$
\left\{
\begin{array}{l}
	\displaystyle{\delta(\theta_6)_{9\, 10}-\eta_7^{9\, 10}=\frac{\partial}{\partial y_{10}}=\tau_{10} -\tau_9 }\\
	\displaystyle{\delta(\theta_6)_{9\, 11}-\eta_7^{9\, 11}=0=\tau_{11} -\tau_9 }
\end{array}
\right.
$$
\hspace{1cm} Hence we have $ \{ \delta(\theta_6)_{9\, 10},\delta(\theta_6)_{9\, 11} \} =\eta_7$.\\
\ \vspace*{3mm}\\
$\circ\  \delta(\theta_7)=\eta_6 -2\eta_7 +\eta_8$
$$
\begin{array}{l}
	\displaystyle{\delta(\theta_7)_{89}=\tilde{\theta_7^9}-\tilde{\theta_7^8}|_{Y_6}=\frac{1+ty_9}{x_9}\frac{\partial}{\partial y_9}-0=\frac{1}{x_9}\frac{\partial}{\partial y_9}=\eta_6^{89} }\\
	\left\{
	\begin{array}{l}
		\displaystyle{\delta(\theta_7)_{9\, 10}=\tilde{\theta_7^{10}}-\tilde{\theta_7^9}|_{Y_7}=-\frac{t+x_{10}}{y_{10}}\frac{\partial}{\partial x_{10}}-\frac{1+ty_9}{x_9}\frac{\partial}{\partial y_9}=-\frac{t+x_{10}}{x_{10}}\frac{\partial}{\partial y_{10}} }\\
		\displaystyle{\delta(\theta_7)_{9\, 11}=\tilde{\theta_7^{11}}-\tilde{\theta_7^9}|_{Y_7}=-\frac{1}{y_{11}}\frac{\partial}{\partial x_{11}}-\frac{1+ty_9}{x_9}\frac{\partial}{\partial y_9}=\frac{t-2x_{11}}{x_{11}(-t+x_{11})}\frac{\partial}{\partial y_{11}} }
	\end{array}
	\right.\\
	\left\{
	\begin{array}{l}
		\displaystyle{\delta(\theta_7)_{11\, 12}=\tilde{\theta_7^{12}}-\tilde{\theta_7^{11}}|_{Y_8}=0-\left(-\frac{1}{y_{11}}\frac{\partial}{\partial x_{11}}\right)=\frac{\partial}{\partial y_{12}}=\eta_{11\, 12} }\\
		\displaystyle{\delta(\theta_7)_{11\, 13}=\tilde{\theta_7^{13}}-\tilde{\theta_7^{11}}|_{Y_8}=0-\left(-\frac{1}{y_{11}}\frac{\partial}{\partial x_{11}}\right)=\frac{1}{x_{13}}\frac{\partial}{\partial y_{13}}=\eta_8^{11\, 13} }
	\end{array}
	\right.
\end{array}
$$
\ \\
\hspace{1cm} Set $\{\tau_9=t\frac{\partial}{\partial x_9},\tau_{10}=-\frac{\partial}{\partial y_{10}},\tau_{11}=0\}\in C^0(N_{Y_7 /S})$.
$$
\left\{
\begin{array}{l}
	\displaystyle{\delta(\theta_7)_{9\, 10}+2\eta_7^{9\, 10}=-\frac{t+x_{10}}{x_{10}}\frac{\partial}{\partial y_{10}}=\tau_{10} -\tau_9}\\
	\displaystyle{\delta(\theta_7)_{9\, 11}+2\eta_7^{9\, 11}=\frac{t}{x_{11}(t-x_{11})}\frac{\partial}{\partial y_{11}}=\tau_{11} -\tau_9}
\end{array}
\right.
$$
\hspace{1cm} Hence we have $ \{ \delta(\theta_7)_{9\, 10},\delta(\theta_7)_{9\, 11} \} =-2\eta_7$.\\
\ \\
$\circ\  \delta(\theta_8)=\eta_7 -2\eta_8$
$$
\begin{array}{l}
	\left\{
	\begin{array}{l}
		\displaystyle{\delta(\theta_8)_{9\, 10}=\tilde{\theta_8^{10}}-\tilde{\theta_8^9}|_{Y_7}=0-0=\eta_7^{9\, 10}}\\
		\displaystyle{\delta(\theta_8)_{9\, 11}=\tilde{\theta_8^{11}}-\tilde{\theta_8^9}|_{Y_7}=\frac{1+a_0y_{11}}{x_{11}}\frac{\partial}{\partial y_{11}}-0=\eta_7^{9\, 11}}
	\end{array}
	\right.\\
	\left\{
	\begin{array}{l}
		\displaystyle{\delta(\theta_8)_{11\, 12}=\tilde{\theta_8^{12}}-\tilde{\theta_8^{11}}|_{Y_8}=-\frac{a_0+x_{12}}{y_{12}}\frac{\partial}{\partial x_{12}}-\frac{1+a_0y_{11}}{x_{11}}\frac{\partial}{\partial y_{11}}=-\frac{a_0+x_{12}}{x_{12}}\frac{\partial}{\partial y_{12}} }\\
		\displaystyle{\delta(\theta_8)_{11\, 13}=\tilde{\theta_8^{13}}-\tilde{\theta_8^{11}}|_{Y_8}=-\frac{1}{y_{13}}\frac{\partial}{\partial x_{13}}-\frac{1+a_0y_{11}}{x_{11}}\frac{\partial}{\partial y_{11}}=\frac{-a_0+2x_{13}}{x_{13}(a_0-x_{13})}\frac{\partial}{\partial y_{13}} }
	\end{array}
	\right.
\end{array}
$$
\ \\
\hspace{1cm} Set $\{\tau_{11}=a_0\frac{\partial}{\partial x_{11}},\tau_{12}=\frac{\partial}{\partial y_{12}},\tau_{13}=0\}\in C^0(N_{Y_8 /S})$.
$$
\left\{
\begin{array}{l}
	\displaystyle{\delta(\theta_8)_{11\, 12}+2\eta_{11\, 12}=\frac{-a_0+x_{12}}{x_{12}}\frac{\partial}{\partial y_{12}}=\tau_{12} -\tau_{11} }\\
	\displaystyle{\delta(\theta_8)_{11\, 13}+2\eta_8^{11\, 13}=\frac{a_0}{x_{13}(a_0-x_{13})}\frac{\partial}{\partial y_{13}}=\tau_{13} -\tau_{11} }
\end{array}
\right.
$$
\hspace{1cm} Hence we have $ \{ \delta(\theta_8)_{11\, 12},\delta(\theta_8)_{11\, 13} \} =-2\eta_8$.\\

Summing up all the computations, we see that the matrix of the linear map $\delta$ is given by 

$$
\left(
\begin{array}{cccccccc}
	-2& 1& 0& 0& 0& 0& 0& 0\\
	1& -2& 1& 0& 0& 0& 0& 0\\
	0& 1& -2& 1& 0& 0& 0& 0\\
	0& 0& 1& -2& 1& 1& 0& 0\\
	0& 0& 0& 1& -2& 0& 0& 0\\
	0& 0& 0& 1& 0& -2& 1& 0\\
	0& 0& 0& 0& 0& 1& -2& 1\\
	0& 0& 0& 0& 0& 0& 1& -2
\end{array}
\right).
$$
\ \\
Since it is well-known that this matrix coincides with the intersection matrix $((Y_i, Y_j))_{1 \leq i,j \leq r}$ of type $\tilde{E_8}$, this completes the proof of Proposition \ref{prop:delta}. \qed
\ \\

\section{Local cohomology of generalized Okamoto--Painlev\'e pair of multiplicative type}\label{sec:mul}
Let $(S,Y)$ be a generalized Okamoto-Painlev\'e pair as in \S \ref{sec:local}.

For a pair $(S,Y)$ of additive type, the result of \S \ref{sec:local} shows the existence of differential equations on $S-Y$ (cf. \cite{STT}). Even for $(S,Y)$ of multiplicative type, if $\dim H^1_D(\Theta_S(-\log D)) \geq 1$, we can derive a differential equation in the same way as in \cite{STT}. Unfortunately, we can prove \linebreak  $H^1_D(\Theta_S(-\log D)) =\{ 0 \}$ for a pair $(S,Y)$ of $\tilde{A_8}$. (For other multiplicative types, we expect that $H^1_D(\Theta_S(-\log D)) =\{ 0 \}$.) This means that, for a pair $(S,Y)$ of multiplicative type, there is no differential equation on $S-Y$ in the sense as in \cite{STT}.

In this section, we will calculate the local cohomology group $H^1_D(\Theta_S(-\log D))$ of pair $(S,Y)$ of $\tilde{A_8}$ type.\\
\ \\
$\bullet$ Construction of $(S,Y)$ of $\tilde{A_8}$ type\\

Now, we consider $(S,Y)$ of $\tilde{A_8}$ type as an example of multiplicative type. \\
According to \cite{Sakai}, any $\tilde{A_8}$-surface is obtained by blowing up $\BP^2$ at the following 9 points given by\\

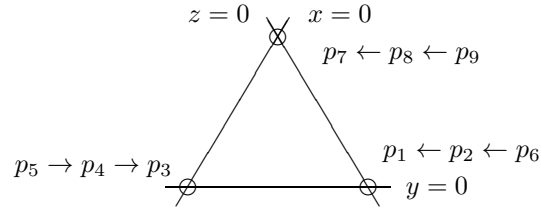
\setlength{\unitlength}{1mm}
\begin{picture}(150,30)(-20,0)
	
	\put(40,5){\line(1,0){30}}
	\put(41.5,2.5){\line(3,5){15}}
	\put(68.5,2.5){\line(-3,5){15}}
	
	\put(43,5){\circle{2}}
	\put(67,5){\circle{2}}
	\put(55,25){\circle{2}}
	
	%%%
	
	\put(59,27){$x=0$}
	\put(72,4){$y=0$}
	\put(43,27){$z=0$}
	
	\put(69,9){$p_1 \gets p_2 \gets p_6$}
	\put(20,7){$p_5 \to p_4 \to p_3$}
	\put(61,22){$p_7 \gets p_8 \gets p_9$}
	
\end{picture}

$$
\begin{array}{l}
	\displaystyle{p_1:(1:0:0) \la p_2:\left( \frac zx,\frac yz \right)=(0,0) \la p_6:\left( \frac zx,\frac{xy}{z^2} \right)=(0,1),}\vspace{3mm}\\
	\displaystyle{p_3:(0:0:1) \la p_4:\left( \frac yz,\frac xy \right)=(0,0) \la p_5:\left( \frac yz,\frac{zx}{y^2} \right)=(0,b),}\vspace{3mm}\\
	\displaystyle{p_7:(0:1:0) \la p_8:\left( \frac xy,\frac zx \right)=(0,0) \la p_9:\left( \frac xy,\frac{yz}{x^2} \right)=(0,c).}
\end{array}
$$

Moreover there exists a $\C^\times$-action on the family of surfaces by
$$
(\mu,(b,c,[x:y:z])) \mapsto (\mu^3 b,\mu^{-3}c,[\mu x:\mu^{-1} y:z]).
$$
By putting $t=bc$, we can normalize this description.\\

We can choose the following coordinate system of $\tilde{A_8}$-surfaces $S_t$ parameterized by t.
\ \\
%%%%%       A8         %%%%%

\setlength{\unitlength}{0,4mm}
\begin{picture}(400,300)(0,0)
	
	\put(200,30){\thicklines\vector(3,1){12}}
	\put(200,30){\thicklines\vector(-3,1){12}}
	
	\put(200,30){\line(3,1){75}}
	\put(200,30){\line(-3,-1){9}}
	
	\put(266,52){\thicklines\vector(-3,-1){12}}
	\put(266,52){\thicklines\vector(1,2){6}}
	
	\put(266,52){\line(1,2){37}}
	\put(267,54){\line(-1,-2){5}}
	
	\put(298,116){\thicklines\vector(-1,-2){6}}
	\put(298,116){\thicklines\vector(-1,4){3}}
	
	\put(298,116){\line(-1,4){17.5}}
	\put(294,132){\line(1,-4){6.5}}
	
	\put(283,176){\thicklines\vector(1,-4){3}}
	\put(283,176){\thicklines\vector(-1,1){9}}
	
	\put(283,176){\line(-1,1){54.5}}
	\put(281,178){\line(1,-1){10}}
	
	\put(237,222){\thicklines\vector(1,-1){9}}
	\put(237,222){\thicklines\vector(-1,0){12}}
	
	\put(237,222){\line(-1,0){85}}
	\put(237,222){\line(1,0){10}}
	
	\put(163,222){\thicklines\vector(1,0){12}}
	\put(163,222){\thicklines\vector(-1,-1){9}}
	
	\put(117,176){\line(1,1){54.5}}
	\put(119,178){\line(-1,-1){10}}
	
	\put(117,176){\thicklines\vector(1,1){9}}
	\put(117,176){\thicklines\vector(-1,-4){3}}
	
	\put(102,116){\line(1,4){17.5}}
	\put(106,132){\line(-1,-4){6.5}}
	
	\put(102,116){\thicklines\vector(1,4){3}}
	\put(102,116){\thicklines\vector(1,-2){6}}
	
	\put(134,52){\line(-1,2){37}}
	\put(133,54){\line(1,-2){5}}
	
	\put(134,52){\thicklines\vector(-1,2){6}}
	\put(134,52){\thicklines\vector(3,-1){12}}
	
	\put(200,30){\line(-3,1){75}}
	\put(200,30){\line(3,-1){9}}	
	
	%%%
	
	\multiput(270,75)(12,-6){6}{\line(2,-1){9}}
	\put(276,72){\thicklines\vector(1,2){6}}
	\put(276,72){\thicklines\vector(2,-1){10}}
	
	\put(332,44){\thicklines\vector(-1,-2){6}}
	\put(332,44){\thicklines\vector(-2,1){10}}

	\multiput(118,99)(-12,-6){6}{\line(-2,-1){9}}
	\put(112,96){\thicklines\vector(1,-2){6}}
	\put(112,96){\thicklines\vector(-2,-1){10}}
	
	\put(56,68){\thicklines\vector(-1,2){6}}
	\put(56,68){\thicklines\vector(2,1){10}}

	\multiput(212,212)(0,13.5){6}{\line(0,1){10}}
	\put(212,222){\thicklines\vector(-1,0){12}}
	\put(212,222){\thicklines\vector(0,1){12}}
	
	\put(212,280){\thicklines\vector(1,0){12}}
	\put(212,280){\thicklines\vector(0,-1){12}}
	
	%%%%%%%%%
	
	\put(268,176){$x_0$}
	\put(273,165){$y_0$}
	\put(284,122){$x_1$}
	\put(282,108){$y_1$}
	\put(260,63){$x_2$}
	\put(250,54){$y_2$}
	\put(205,38){$x_3$}
	\put(186,38){$y_3$}
	\put(145,54){$x_4$}
	\put(133,63){$y_4$}
	\put(110,108){$x_5$}
	\put(108,122){$y_5$}
	\put(118,165){$x_6$}
	\put(125,176){$y_6$}
	\put(158,208){$x_7$}
	\put(172,213){$y_7$}
	\put(224,213){$x_8$}
	\put(236,208){$y_8$}
	\put(270,84){$x_9$}
	\put(288,70){$y_9$}
	\put(120,86){$x_{10}$}
	\put(92,95){$y_{10}$}
	\put(195,213){$x_{11}$}
	\put(215,235){$y_{11}$}
	\put(320,53){$x_{12}$}
	\put(330,30){$y_{12}$}
	\put(65,64){$x_{13}$}
	\put(36,80){$y_{13}$}
	\put(197,268){$x_{14}$}
	\put(222,284){$y_{14}$}
	
	%%%%%%%%%%%
	
	\put(295,145){$Y_1$}
	\put(256,88){$Y_2$}
	\put(230,30){$Y_3$}
	\put(163,30){$Y_4$}
	\put(136,88){$Y_5$}
	\put(95,145){$Y_6$}
	\put(126,200){$Y_7$}
	\put(195,198){$Y_8$}
	\put(265,200){$Y_9$}
	\put(195,250){$E_3$}
	\put(70,85){$E_2$}
	\put(292,50){$E_1$}
	\put(270,250){$E_i \simeq \BP^1,\ E_i^2=-1,$}
	\put(270,237){$Y_i \simeq \BP^1,\ Y_i^2=-2.$}
\end{picture}

\begin{eqnarray}
	U_i & = & \Spec \C \left[x_i,y_i\right] \ \cong \ {\C}^2\ \ (i=0,1,...,8).\nonumber\\
	U_9 & = & \Spec \C \left[x_9,y_9,\frac{1}{t+x_9}\right] \ \cong \ {\C}^2-\{t+x_9=0\}.\nonumber\\
	U_i& = & \Spec \C \left[x_i,y_i,\frac{1}{1+x_i}\right] \ \cong \ {\C}^2-\{1+x_i=0\}\ \ (i=10,11).\nonumber\\
	U_{12} & = & \Spec \C \left[u_{12},v_{12},\frac{1}{t+u_{12} v_{12}}\right] \ \cong \ {\C}^2-\{ t+u_{12} v_{12}=0\}.\nonumber\\
	U_i& = & \Spec \C \left[u_i,v_i,\frac{1}{1+u_i v_i}\right] \ \cong \ {\C}^2-\{ 1+u_i v_i=0\}\ \ (i=13,14).\nonumber
\end{eqnarray}
where $t \in \C ^\times$.

$$
\begin{array}{ll}
	Y_1=\{x_0=0,y_1=0 \},& Y_2=\{x_1=0,y_2=0,y_9=0 \},\vspace{2mm}\\
	Y_3=\{x_2=0,y_3=0 \},& Y_4=\{x_3=0,y_4=0 \},\vspace{2mm}\\
	Y_5=\{x_4=0,y_5=0,y_{10}=0 \},& Y_6=\{x_5=0,y_6=0 \},\vspace{2mm}\\
	Y_7=\{x_6=0,y_7=0\},& Y_8=\{x_7=0,y_8=0,y_{11}=0 \},\vspace{2mm}\\
	Y_9=\{x_8=0,y_0=0\},&
\end{array}
$$

$$ S_t=\bigcup_{i=0}^{14}U_i,$$
$$Y_t=\sum_{i=1}^{9}Y_i,\ \ \ D=Y,$$
$$S_t-Y_t=U_{12} \cup U_{13} \cup U_{14}.$$

For later use, we need only the coordinate transformations near each component $Y_i$.
Here we will list up the coordinate transformations only for a neighborhood of each $Y_i$.
$$
\hspace*{-2mm}
\begin{array}{lll}
Y_1:\left\{
\begin{array}{l}
	x_0={x_1}^2 y_1\\
	y_0=\displaystyle{\frac1{x_1}}
\end{array}
\right.\hspace{8mm}
&
Y_2:
\left\{
\begin{array}{l}
	x_1=x_2 y_2\\
	y_1=\displaystyle{\frac1{x_2}}
\end{array}
\right.
\left\{
\begin{array}{l}
	x_1= x_9 (t+x_9) y_9\\
	\displaystyle{y_1=\frac1{(t+x_9)}}
\end{array}
\right.\hspace{8mm}
&
Y_3:
\left\{
\begin{array}{l}
	x_2={x_3}^2 y_3\\
	y_2=\displaystyle{\frac1{x_3}}
\end{array}
\right.
\vspace{3mm}\\
Y_4:
\left\{
\begin{array}{l}
	x_3={x_4}^2 y_4\\
	y_3=\displaystyle{\frac1{x_4}}
\end{array}
\right.\hspace{8mm}
&
Y_5:
\left\{
\begin{array}{l}
	x_4=x_5 y_5\\
	y_4=\displaystyle{\frac1{x_5}}
\end{array}
\right.
\left\{
\begin{array}{l}
	x_4=x_{10} (1+ x_{10}) y_{10}\\
	y_4=\displaystyle{\frac1{(1+ x_{10})}}
\end{array}
\right.\hspace{8mm}
&
Y_6:
\left\{
\begin{array}{l}
	x_5={x_6}^2 y_6\\
	y_5=\displaystyle{\frac1{x_6}}
\end{array}
\right.
\vspace{3mm}\\
Y_7:
\left\{
\begin{array}{l}
	x_6={x_7}^2 y_7\\
	y_6=\displaystyle{\frac1{x_7}}
\end{array}
\right.\hspace{8mm}
&
Y_8:
\left\{
\begin{array}{l}
	x_7=x_8 y_8\\
	y_7=\displaystyle{\frac1{x_8}}
\end{array}
\right.
\left\{
\begin{array}{l}
	x_7= x_{11} (1+x_{11}) y_{11}\\
	y_7=\displaystyle{\frac1{(1+ x_{11})}}
\end{array}
\right.\hspace{8mm}
&
Y_9:
\left\{
\begin{array}{l}
	x_8={x_0}^2 y_0\\
	y_8=\displaystyle{\frac1{x_0}}.
\end{array}
\right.
\end{array}
$$

This gives a generalized Okamoto--Painlev\'e pair $(S_t,Y_t)$ of type ${\tilde A}_8$.

\begin{Proposition}\label{prop:multi}
Let $(S_t, Y_t)$ be as above and set $D_t=(Y_t)_{red}$.\\
If \ $-t \in \C$ is not a root of unity, then we have
$$
H^1_{D_t}(\Theta_{S_t}(-\log D_t))=0.
$$
\end{Proposition}
\begin{Remark}
{\rm We expect that, if $\, -t$ is a root of unity then $(S_t,Y_t)$ is of fibered type.}
\end{Remark}

{\it Proof.}
From the diagram
$$
\begin{array}{ccccccc}
	&&&0&&0&\\
	&&& \downarrow &&\downarrow&\\
	0 \lra & \Theta_S(-\log D) & \lra & \Theta_S(-\log D) \otimes \cO_S((n-1)D) & \lra & \Theta_S(-\log D) \otimes N_{(n-1)D} & \lra 0\\
	& \parallel && \downarrow  && \downarrow \mu &\\
	0 \lra & \Theta_S(-\log D) & \lra & \Theta_S(-\log D) \otimes \cO_S(nD) & \lra & \Theta_S(-\log D) \otimes N_{nD} & \lra 0\\
	&&& \downarrow && \downarrow &\\
	&&& \Theta_S(-\log D) \otimes N^{\otimes n}_D & \simeq & \coker \mu &\\
	&&& \downarrow && \downarrow &\\
	&&&0&&0&,\\
\end{array}
$$
we obtain the exact sequence
$$
0 \lra \Theta_S(-\log D) \otimes N_{(n-1)D} \lra \Theta_S(-\log D) \otimes N_{nD} \lra \Theta_S(-\log D) \otimes N^{\otimes n}_D \lra 0.
$$
Therefore, we have the following sequence for each $n \geq 2$. 
$$
0 \lra H^0(\Theta_S(-\log D) \otimes N_{(n-1)D}) \lra H^0(\Theta_S(-\log D) \otimes N_{nD}) \lra H^0(\Theta_S(-\log D) \otimes N^{\otimes n}_D) .
$$
If $H^0(\Theta_S(-\log D) \otimes N_D^{\otimes n})=0$ for any $n \geq 1$,
then we have
$$
H^0(\Theta_S(-\log D) \otimes N_{(n-1)D}) \simeq H^0(\Theta_S(-\log D) \otimes N_{nD}).
$$
Since $H^0(\Theta_S(-\log D) \otimes N_D)=0$, we have 
$$
H^0(\Theta_S(-\log D) \otimes N_{nD})=0\ \ (n \geq 1).
$$
Therefore, noting that (cf. \S \ref{sec:local})
$$
H^1_D(S, \Theta_S( - \log D)) = \varinjlim H^0( \Theta_S( - \log D) \otimes N_{nD}),
$$
we obtain
$$
H^1_D(S, \Theta_S( - \log D)) = 0.
$$\vspace{2mm}

Now let us calculate $H^0(\Theta_S(-\log D) \otimes N_D^{\otimes n})$.

Recall that we have the exact sequence
$$
0 \lra \nu_\ast(\oplus_{i=1}^9 N_{Y_i} ) \lra \Theta_S( - \log D) \otimes N_D  \lra \Theta_D \otimes N_D \lra 0,
$$
where $\nu : \tilde{D}=\coprod^9_{i=1} Y_i \to D$ is the normalization map (Lemma \ref{lem:kuma} (\ref{eq:thetalog})).
Tensoring this sequence with $N_D^{\otimes n-1}$, we obtain
$$
0 \lra \nu_\ast(\oplus_{i=1}^9 N_{Y_i}) \otimes N_D^{\otimes n-1} \lra \Theta_S( - \log D) \otimes N_D^{\otimes n}  \lra \Theta_D \otimes N_D^{\otimes n} \lra 0.
$$
Therefore we have
$$
\begin{array}{l}
H^0(D,\nu_\ast(\oplus_{i=1}^9 N_{Y_i}) \otimes N_D^{\otimes n-1}) \ra H^0(D,\Theta_S( - \log D) \otimes N_D^{\otimes n}) \ra H^0(D,\Theta_D \otimes N_D^{\otimes n} ) \ra \\
H^1(D,\nu_\ast(\oplus_{i=1}^9 N_{Y_i}) \otimes N_D^{\otimes n-1}).
\end{array}
$$

In $\tilde{A_8}$ case, we see that $t_i = 2$, where $t_i$ is the number of intersections of $Y_i$ with other components. From this, we obtain
$$
\Theta_D \simeq \nu_\ast(\oplus^9_{i=1} \Theta_{Y_i} ( -t_i)) \simeq  \nu_\ast(\oplus^9_{i=1} \cO_{Y_i} (2 - t_i)) \simeq \nu_\ast(\oplus^9_{i=1} \cO_{Y_i}) .
$$
Since $\nu$ is a finite morphism, we have
$$
\begin{array}{l}
H^i(D,\nu_\ast(\oplus_{i=1}^9 N_{Y_i}) \otimes N_D^{\otimes n-1}) \simeq H^i(\tilde{D},(\oplus_{i=1}^9 N_{Y_i}) \otimes \nu^\ast(N_D^{\otimes n-1}))\ \ \ (i=0,1),\vspace{2mm}\\
H^0(D,\Theta_D \otimes N_D^{\otimes n}) \simeq H^0(D,\nu_\ast(\oplus^9_{i=1} \cO_{Y_i}) \otimes N_D^{\otimes n} ) \simeq H^0(\tilde{D},(\oplus^9_{i=1} \cO_{Y_i}) \otimes \nu^\ast(N_D^{\otimes n})).
\end{array}
$$
Note that $D \cdot Y_i = 2-t_i=0$, we see that ${N_D}|_{Y_i} \simeq \cO_{Y_i}$. Therefore we have
$$
(\oplus_{i=1}^9 N_{Y_i}) \otimes \nu^\ast(N_D^{\otimes n-1}) \simeq \oplus_{i=1}^9 N_{Y_i},\ \ \  (\oplus^9_{i=1} \cO_{Y_i}) \otimes \nu^\ast(N_D^{\otimes n}) \simeq \oplus^9_{i=1} \cO_{Y_i}.
$$

Summing up the arguments above, we obtain
$$
\begin{array}{ccccccc}
	0 & \lra & H^0(\Theta_S( - \log D) \otimes N_D^{\otimes n}) & \lra & \oplus_{i=1}^9 H^0(Y_i, \cO_{Y_i}) & \stackrel{\delta}{\lra} &  \oplus_{i=1}^9 H^1(N_{Y_i})\\
	  &  &  &  & \parallel &  & \parallel \\
	  &  &  &  & \C^9 & \stackrel{\delta}{\lra} & \C^9.
\end{array} 
$$

We will analyse the edge homomorphism $\delta$ by using the \v{C}ech cocycles. \vspace{2mm}

Noting  that 
\begin{equation}\label{eq:tori}
H^0(\Theta_D \otimes N_D^{\otimes n}) \simeq H^0(\Theta_{Y_i}(-2) \otimes N_D^{\otimes n}) \simeq \C^9.
\end{equation}
For each $i\ (1 \le i \le 9)$, we introduce a generator $\theta_i$ of the cohomology group in (\ref{eq:tori}) as follows.\\

\begin{center}
\begin{tabular}{|c|l|}\hline
&\vspace{-2mm}\\
$\theta_1$ & 
$\left\{\theta_1^0=\frac{y_0}{x_0^n y_0^n}\frac{\partial}{\partial y_0}\ \mbox{on} \ U_0\cap Y_1 ,\ \ 
\theta_1^1=-\frac{x_1}{x_1^n y_1^n}\frac{\partial}{\partial x_1}\ \mbox{on} \ U_1\cap Y_1 \right\}$\\
&\vspace{-2mm}\\
\hline
&\vspace{-2mm}\\
$\theta_2$ & 
$\left\{
\begin{array}{c}
\theta_2^1=\frac{(1-ty_1)^n y_1}{x_1^n y_1^n}\frac{\partial}{\partial y_1} \ \mbox{on} \ U_1\cap Y_2,\ \ 
\theta_2^2=-\frac{(x_2-t)^n x_2}{x_2^n y_2^n}\frac{\partial}{\partial x_2} \ \mbox{on} \ U_2\cap Y_2,\vspace{1mm}\\
\theta_2^9=-\frac{1}{(t+ x_9)^{n-1} y_9^n}\frac{\partial}{\partial x_9}\ \mbox{on} \ U_9\cap Y_2 
\end{array}
\right\}$\\
&\vspace{-2mm}\\
\hline
&\vspace{-2mm}\\
$\theta_3$ & 
$\left\{\theta_3^2=\frac{y_2}{x_2^n y_2^n}\frac{\partial}{\partial y_2}\ \mbox{on} \ U_2\cap Y_3,\ \ 
\theta_3^3=-\frac{x_3}{x_3^n y_3^n}\frac{\partial}{\partial x_3}\ \mbox{on} \ U_3\cap Y_3\right\}$\\
&\vspace{-2mm}\\
\hline
&\vspace{-2mm}\\
$\theta_4$ & 
$\left\{\theta_4^3=\frac{y_3}{x_3^n y_3^n}\frac{\partial}{\partial y_3}\ \mbox{on} \ U_3\cap Y_4,\ \ 
\theta_4^4=-\frac{x_4}{x_4^n y_4^n}\frac{\partial}{\partial x_4}\ \mbox{on} \ U_4\cap Y_4\right\}$\\
&\vspace{-2mm}\\
\hline
&\vspace{-2mm}\\
$\theta_5$ & 
$\left\{
\begin{array}{c}
\theta_5^4=\frac{(1-x_4)^n y_4}{x_4^n y_4^n}\frac{\partial}{\partial y_4}\ \mbox{on} \ U_4\cap Y_5 ,\ \ 
\theta_5^5=-\frac{(x_5-1)^n x_5}{x_5^n y_5^n}\frac{\partial}{\partial x_5}\ \mbox{on} \ U_5\cap Y_5,\vspace{1mm}\\
\theta_5^{10}=-\frac{1}{(1+ x_{10})^{n-1} y_{10}^n}\frac{\partial}{\partial x_{10}}\ \mbox{on} \ U_{10}\cap Y_5
\end{array}
\right\}$\\
&\vspace{-2mm}\\
\hline
&\vspace{-2mm}\\
$\theta_6$ & 
$\left\{\theta_6^5=(-1)^n \frac{y_5}{x_5^n y_5^n}\frac{\partial}{\partial y_5}\ \mbox{on} \ U_5\cap Y_6,\ \ 
\theta_6^6=-(-1)^n\frac{x_6}{x_6^n y_6^n}\frac{\partial}{\partial x_6}\ \mbox{on} \ U_6\cap Y_6\right\}$\\
&\vspace{-2mm}\\
\hline
&\vspace{-2mm}\\
$\theta_7$ & 
$\left\{\theta_7^6=(-1)^n \frac{y_6}{x_6^n y_6^n}\frac{\partial}{\partial y_6}\ \mbox{on} \ U_6\cap Y_7,\ \ 
\theta_7^7=-(-1)^n \frac{x_7}{x_7^n y_7^n}\frac{\partial}{\partial x_7}\ \mbox{on} \ U_7\cap Y_7\right\}$\\
&\vspace{-2mm}\\
\hline
&\vspace{-2mm}\\
$\theta_8$ & 
$\left\{
\begin{array}{c}
\theta_8^7=(-1)^n \frac{(1-x_7)^n y_7}{x_7^n y_7^n}\frac{\partial}{\partial y_7}\ \mbox{on} \ U_7\cap Y_8,\ \ 
\theta_8^8=-(-1)^n \frac{(x_8-1)^n x_8}{x_8^n y_8^n}\frac{\partial}{\partial x_8}\ \mbox{on} \ U_8\cap Y_8,\vspace{1mm}\\
\theta_8^{11}=-(-1)^n \frac{1}{(1+ x_{11})^{n-1} y_{11}^n}\frac{\partial}{\partial x_{11}}\ \mbox{on} \ U_{11}\cap Y_8 
\end{array}
\right\}$\\
&\vspace{-2mm}\\
\hline
&\vspace{-2mm}\\
$\theta_9$ & 
$\left\{\theta_9^8=\frac{y_8}{x_8^n y_8^n}\frac{\partial}{\partial y_8}\ \mbox{on} \ U_8\cap Y_9,\ \ 
\theta_9^0=-\frac{x_0}{x_0^n y_0^n}\frac{\partial}{\partial x_0}\ \mbox{on} \ U_0\cap Y_9\right\}$\vspace{-2mm}\\
&\\
\hline
\end{tabular}
\end{center}

On the other hand, for each $i\ (1 \le i \le 9)$, we have a generator $\eta_i \in H^1(Y_i,N_{Y_i/S} \otimes N_D^{\otimes n-1})$ as follows.

\begin{center}
\begin{tabular}{|c|l|}\hline
&\vspace{-2mm}\\
$\eta_1$ & 
$\left\{\eta_1^{01}=\frac{1}{x_1x_1^{n-1} y_1^{n-1}}\frac{\partial}{\partial y_1}\ \mbox{on} \ U_0\cap \ U_1\cap Y_1\right\}$\\
&\vspace{-2mm}\\
\hline
&\vspace{-2mm}\\
$\eta_2$ & 
$\left\{\eta_2^{12}=\frac{1}{y_2^{n-1}}\frac{\partial}{\partial y_2}\ \mbox{on} \ U_1\cap \ U_2\cap Y_2,\ \eta_2^{29}=0\ \mbox{on} \ U_2\cap \ U_9\cap Y_2\right\}$\\
&\vspace{-2mm}\\
\hline
&\vspace{-2mm}\\
$\eta_3$ & 
$\left\{\eta_3^{23}=\frac{1}{x_3x_3^{n-1} y_3^{n-1}}\frac{\partial}{\partial y_3}\ \mbox{on} \ U_2\cap \ U_3\cap Y_3\right\}$\\
&\vspace{-2mm}\\
\hline
&\vspace{-2mm}\\
$\eta_4$ & 
$\left\{\eta_4^{34}=\frac{1}{x_4x_4^{n-1} y_4^{n-1}}\frac{\partial}{\partial y_4}\ \mbox{on} \ U_3\cap \ U_4\cap Y_4\right\}$\\
&\vspace{-2mm}\\
\hline
&\vspace{-2mm}\\
$\eta_5$ & 
$\left\{\eta_5^{45}=\frac{1}{y_5^{n-1}}\frac{\partial}{\partial y_5}\ \mbox{on} \ U_4\cap \ U_5\cap Y_5,\ \eta_5^{5\, 10}=0\ \mbox{on} \ U_5\cap \ U_{10}\cap Y_5\right\}$\\
&\vspace{-2mm}\\
\hline
&\vspace{-2mm}\\
$\eta_6$ & 
$\left\{\eta_6^{56}=(-1)^n \frac{1}{x_6x_6^{n-1} y_6^{n-1}}\frac{\partial}{\partial y_6}\ \mbox{on} \ U_5\cap \ U_6\cap Y_6\right\}$\\
&\vspace{-2mm}\\
\hline
&\vspace{-2mm}\\
$\eta_7$ & 
$\left\{\eta_7^{67}=(-1)^n \frac{1}{x_7x_7^{n-1} y_7^{n-1}}\frac{\partial}{\partial y_7}\ \mbox{on} \ U_6\cap \ U_7\cap Y_7\right\}$\\
&\vspace{-2mm}\\
\hline
&\vspace{-2mm}\\
$\eta_8$ & 
$\left\{\eta_8^{78}=(-1)^n \frac{1}{y_8^{n-1}}\frac{\partial}{\partial y_8}\ \mbox{on} \ U_7\cap \ U_8\cap Y_8,\ \eta_8^{8\, 11}=0\ \mbox{on} \ U_8\cap \ U_{11}\cap Y_8\right\}$\\
&\vspace{-2mm}\\
\hline
&\vspace{-2mm}\\
$\eta_9$ & 
$\left\{\eta_9^{80}=\frac{1}{x_0x_0^{n-1} y_0^{n-1}}\frac{\partial}{\partial y_0}\ \mbox{on} \ U_8\cap \ U_0\cap Y_9\right\}$\vspace{-2mm}\\
&\\
\hline
\end{tabular}
\end{center}

\ \\
We take $\{\theta_i\}$ and $\{\eta_i\}$ as basis of $H^0(\Theta_D \otimes N_D^{\otimes n})$ and $\oplus_{i=1}^9 H^1(N_{Y_i /S} \otimes N_D^{\otimes n-1})$ respectively. By using these bases, we compute the matrix representing the connecting homomorphisn $\delta$.

For that purpose, let us lift 0-cocycle $\theta_1$ to 0-cochains of $\Theta_S(- \log D) \otimes N_D$ as 
$$
\tilde{\theta_1^0}=\frac{y_0}{x_0^n y_0^n}\frac{\partial}{\partial y_0}
  \ {\rm on}\ U_0,\ \ \ 
\tilde{\theta_1^1}=-\frac{x_1}{x_1^n y_1^n}\frac{\partial}{\partial x_1}  \ {\rm on}\ U_1,
$$
$$\tilde{\theta_1^i}=0\ {\rm on}\ U_i\ (i=2,3,\cdots,14).$$
Other 0-cocycles can be lifted in a similar way.\\

We first compute $\delta(\theta_1)$.\\
\ \\
From the definition of $\delta$, we have $\delta(\theta_1)=\{ \delta(\theta_1)_{ij} \ {\rm on}\ U_i \cap U_j \cap D \}$ with
$$
\begin{array}{l}
	\displaystyle{\delta(\theta_1)_{80}=-\tilde{\theta_1^8}+\tilde{\theta_1^0}|_{Y_9}=-0+\frac{y_0}{x_0^n y_0^n}\frac{\partial}{\partial y_0}=\eta_9^{80}}\\
	\displaystyle{\delta(\theta_1)_{01}=-\tilde{\theta_1^0}+\tilde{\theta_1^1}|_{Y_1}=-\frac{y_0}{x_0^n y_0^n}\frac{\partial}{\partial y_0}+\left(-\frac{x_1}{x_1^n y_1^n}\frac{\partial}{\partial x_1}\right)=-2\frac{1}{x_1x_1^{n-1} y_1^{n-1}}\frac{\partial}{\partial y_1}=-2\eta_1^{01}}\\
	\left\{
	\begin{array}{l}
		\displaystyle{\delta(\theta_1)_{12}=-\tilde{\theta_1^1}+\tilde{\theta_1^2}|_{Y_2}=-\left(-\frac{x_1}{x_1^n y_1^n}\frac{\partial}{\partial x_1}\right)+0=\frac{1}{y_2^{n-1}}\frac{\partial}{\partial y_2}=\eta_2^{12}}\\
		\delta(\theta_1)_{29}=-\tilde{\theta_1^2}+\tilde{\theta_1^9}|_{Y_2}=-0+0=\eta_2^{29}
	\end{array}
	\right.
\end{array}
$$
Other $\delta(\theta_1)_{ij}$'s are zero.\\
\ \\
Obviously $\delta(\theta_1)=\eta_9-2\eta_1+\eta_2$.

Other $\delta(\theta_i)$'s can be treated in a similar way. In what follows, we just list up a few results of computations.\\

$\circ\  \delta(\theta_2)=\eta_1-2\eta_2 + (-t)^n \eta_3$
$$
\begin{array}{l}
	\displaystyle{\delta(\theta_2)_{01}=-\tilde{\theta_2^0}+\tilde{\theta_2^1}|_{Y_1}=-0+\frac{(1-t y_1)^n y_1}{x_1^n y_1^n}\frac{\partial}{\partial y_1}=\frac{1}{x_1x_1^{n-1} y_1^{n-1}}\frac{\partial}{\partial y_1}=\eta_1^{01}}\ (y_1=0\ {\rm on}\ U_1 \cap Y_1)\\
	\hspace{-2mm}\left\{
	\begin{array}{l}
		\hspace{-1.5mm}\displaystyle{\delta(\theta_2)_{12}=-\tilde{\theta_2^1}+\tilde{\theta_2^2}|_{Y_2}=-\frac{(1-t y_1)^n y_1}{x_1^n y_1^n}\frac{\partial}{\partial y_1}+\left(-\frac{(x_2-t)^n x_2}{x_2^n y_2^n}\frac{\partial}{\partial x_2}\right)=-\frac{(x_2-t)^n}{x_2x_2^{n-1} y_2^{n-1}}\frac{\partial}{\partial y_2}}\\
		\hspace{-1.5mm}\displaystyle{\delta(\theta_2)_{29}=-\tilde{\theta_2^2}+\tilde{\theta_2^9}|_{Y_2}=\hspace{-0.5mm}-\hspace{-1mm}\left(\hspace{-1mm}-\frac{(x_2-t)^n x_2}{x_2^n y_2^n}\frac{\partial}{\partial x_2}\hspace{-0.5mm}\right)\hspace{-1mm}+\hspace{-1mm}\left(\hspace{-1mm}-\frac{1}{(t+ x_9)^{n-1} y_9^n}\frac{\partial}{\partial x_9}\hspace{-0.5mm}\right)\hspace{-1mm}=\hspace{-0.5mm}-\frac{1}{x_9 (t+x_9)^{n-1} y_9^{n-1}}\frac{\partial}{\partial y_9}}
	\end{array}
	\right.\\

	\displaystyle{\delta(\theta_2)_{23}=-\tilde{\theta_2^2}+\tilde{\theta_2^3}|_{Y_3}=-\left(-\frac{(x_2-t)^n x_2}{x_2^n y_2^n}\frac{\partial}{\partial x_2}\right)\frac{\partial}{\partial x_2}+0=(-t)^n \frac{1}{x_3x_3^{n-1} y_3^{n-1}}\frac{\partial}{\partial y_3}=(-t)^n \eta_3^{23}}
\end{array}
$$
\ \\
Set $\{\tau_1=\frac{(1-ty_1)^n+(1-ty_1)^{n-1}-2}{y_1x_1^{n-1} y_1^{n-1}}\frac{\partial}{\partial x_1},\tau_2=\frac{(x_2-t)^{n-1}}{x_2^{n-1} y_2^{n-1}}\frac{\partial}{\partial y_2},\tau_9=0\}\in C^0(N_{Y_2/M}\otimes N_D^{\otimes n-1})$.
$$
\left\{
\begin{array}{l}
\displaystyle{\delta(\theta_2)_{12}+2\eta_2^{12}=-\frac{(x_2-t)^n}{x_2x_2^{n-1} y_2^{n-1}}\frac{\partial}{\partial y_2}+2\frac{1}{y_2^{n-1}}\frac{\partial}{\partial y_2}=-\tau_1+\tau_2}\\

\displaystyle{\delta(\theta_2)_{29}+2\eta_2^{29}=-\frac{1}{x_9 (t+x_9)^{n-1} y_9^{n-1}}\frac{\partial}{\partial y_9}-0=-\tau_2+\tau_9}
\end{array}
\right.
$$

Thus we have $ \{ \delta(\theta_2)_{12},\delta(\theta_2)_{29} \} =-2\eta_2$.\\
\ \\

$\circ\  \delta(\theta_3)=\frac1{(-t)^n} \eta_2 -2\eta_3+\eta_4$
$$
\begin{array}{l}
	\left\{
	\begin{array}{l}
		\displaystyle{\delta(\theta_3)_{12}=-\tilde{\theta_3^1}+\tilde{\theta_3^2}|_{Y_2}=-0+\frac{y_2}{x_2^n y_2^n}\frac{\partial}{\partial y_2}=\frac{y_2}{x_2^n y_2^n}\frac{\partial}{\partial y_2}}\\
		\displaystyle{\delta(\theta_3)_{29}=-\tilde{\theta_3^2}+\tilde{\theta_3^9}|_{Y_2}=-\frac{y_2}{x_2^n y_2^n}\frac{\partial}{\partial y_2}+0=-\frac{1}{x_9^n (t+x_9)^n y_9^{n-1}}\frac{\partial}{\partial y_9}}
	\end{array}
	\right.\\

	\displaystyle{\delta(\theta_3)_{23}=-\tilde{\theta_3^2}+\tilde{\theta_3^3}|_{Y_3}=-\frac{y_2}{x_2^n y_2^n}\frac{\partial}{\partial y_2}+\left(-\frac{x_3}{x_3^n y_3^n}\frac{\partial}{\partial x_3}\right)=-2\frac{1}{x_3x_3^{n-1} y_3^{n-1}}\frac{\partial}{\partial y_3}=-2\eta_3^{23}}\\
	\displaystyle{\delta(\theta_3)_{34}=-\tilde{\theta_3^3}+\tilde{\theta_3^4}|_{Y_4}=-\left(-\frac{x_3}{x_3^n y_3^n}\frac{\partial}{\partial x_3}\right)+0=\frac{1}{x_4x_4^{n-1} y_4^{n-1}}\frac{\partial}{\partial y_4}=\eta_4^{34}}
\end{array}
$$
\ \\
Set $\{\tau_1=\frac{1-(1-ty_1)^n}{(-t)^n y_1x_1^{n-1} y_1^{n-1}}\frac{\partial}{\partial x_1},\tau_2=\frac{(-t)^n-(x_2-t)^n}{(-t)^n x_2x_2^{n-1} y_2^{n-1}}\frac{\partial}{\partial y_2},\tau_9=-\frac{1}{(-t)^n (x_9+t)^n y_9^{n-1}}\frac{\partial}{\partial y_9}\}\in C^0(N_{Y_2/M}\otimes N_D^{\otimes n-1})$.
$$
\left\{
\begin{array}{l}
\displaystyle{\delta(\theta_3)_{12}-\frac1{(-t)^n}\eta_2^{12}=\frac{y_2}{x_2^n y_2^n}\frac{\partial}{\partial y_2}}-\frac1{(-t)^n}\frac{1}{y_2^{n-1}}\frac{\partial}{\partial y_2}=-\tau_1+\tau_2 \\

\displaystyle{\delta(\theta_3)_{29}-\frac1{(-t)^n}\eta_2^{29}=-\frac{1}{x_9^n (t+x_9)^n y_9^{n-1}}\frac{\partial}{\partial y_9}}-0=-\tau_2+\tau_9
\end{array}
\right.
$$

Thus we have $ \{ \delta(\theta_3)_{12},\delta(\theta_3)_{29} \} =\frac1{(-t)^n}\eta_2$.\\
\ \\
\ \\
Summing up all the computations, we see that the matrix of the linear map $\delta$ is given by

$$
\delta = 
\left(
\begin{array}{ccccccccc}
	-2& 1& 0& 0& 0& 0& 0& 0& 1\\
	1& -2& (-t)^{-n} & 0& 0& 0& 0& 0& 0\\
	0& (-t)^n & -2& 1& 0& 0& 0& 0& 0\\
	0& 0& 1& -2& 1& 0& 0& 0& 0\\
	0& 0& 0& 1& -2& 1& 0& 0& 0\\
	0& 0& 0& 0& 1& -2& 1& 0& 0\\
	0& 0& 0& 0& 0& 1& -2& 1& 0\\
	0& 0& 0& 0& 0& 0& 1& -2& 1\\
	1& 0& 0& 0& 0& 0& 0& 1& -2
\end{array}
\right)
$$
Since
$$ 
\det \delta =\frac{((-t)^n-1)^2}{(-t)^n} ,
$$
we see
$$
\rank \delta =
\left\{\begin{array}{l}
	9\ \ \ (-t {\rm \ is\ not\ a\ root\ of\ unity})\\
	8\ \ \ (-t {\rm \ is\ a\ root\ of\ unity}).
\end{array}\right.
$$

Therefore, if \ $-t$ is not a root of unity, then
$$
H^0(\Theta_S(-\log D) \otimes N_D^{\otimes n})=0.
$$
\qed \vspace{5mm}

{\bf Acknowledgements}

The author is deeply grateful to Professor Masa-Hiko Saito, who gave support to the whole of this work. We would also like to thank Kenji Iohara for many useful advices.

\end{document}